\documentclass[11pt]{amsart}
\usepackage{fullpage}
\usepackage{amssymb}
\usepackage{amsmath}
\usepackage{amsxtra}
\usepackage{amscd}
\usepackage{graphicx}
\usepackage{amsfonts}
\usepackage{pb-diagram}
\usepackage{float}
\usepackage{enumerate}

\numberwithin{equation}{section}

\newtheorem{thm}{Theorem}

\newtheorem{cor}{Corollary}
\newtheorem{prop}{Proposition}

\newtheorem{defn}{Definition}

\newtheorem{rem}{Remark}
\newtheorem{note}{Note}

\begin{document}

\title{On the framization of knot algebras}

\author{Jes\'us Juyumaya}
\address{J. Juyumaya: Instituto de Matem\'aticas,  Universidad de Valpara\'{\i}so \\
Gran Breta\~na 1091, Valpara\'{\i}so, Chile.}
\email{juyumaya@uvach.cl}

\author{Sofia Lambropoulou}
\address{S. Lambropoulou: Department of Mathematics,
National Technical University of Athens,
Zografou campus, GR--157 80 Athens, Greece.}
\email{sofia@math.ntua.gr}
\urladdr{http://www.math.ntua.gr/$\sim$sofia}

\thanks{This research has been co--financed by the European Union (European
Social Fund -- ESF) and Greek national funds through the Operational Program
``Education and Lifelong Learning" of the National Strategic Reference
Framework (NSRF) -- Research Funding Program: THALIS. Moreover, the first author was partially supported by Fondecyt 1141254  and Dipuv No. 1/2011. }

\keywords{framed links, framed braids, Yokonuma--Hecke algebra, Markov trace, knot invariants, singular links, Temperley--Lieb algebra, BMW algebra, $B$--type Hecke algebras, affine Hecke algebra, singular Hecke algebra}

\subjclass[2010]{57M27, 57M25, 20F38, 20F36, 20C08}

\date{}

\begin{abstract}
This paper presents results on the framization of some knot algebras, defined by the authors.  We explain the motivations of the concept of framization, coming from the Yokonuma--Hecke algebras, as well as recent results on the framization of the Temperley--Lieb algebra.
Finally,  we propose framizations for other knot algebras such as the BMW algebra, the $B$--type related Hecke algebras and the singular Hecke algebra.
\end{abstract}

\maketitle

\section*{Introduction}

Modular framization (or simply framization) is a mechanism proposed recently by the authors and it consists in constructing a non--trivial extension of a knot algebra via the addition of framing generators. In this way we obtain a new algebra which is related to framed braids and framed knots. 

\smallbreak

By knot algebra we mean an algebra that is applied to the construction of invariants of knots and links. We are focused  on the framization of those knot algebras that define invariants of knots via the Jones' construction. More precisely, such a knot algebra $A$ is   a triplet $(A, \pi, \tau)$, where $\pi$ is a representation of a braid group in $A$ and $\tau$ a Markov trace defined on $A$. The invariant obtained by  the knot algebra $A$ is constructed essentially from the composition $\tau\circ\pi$ after  re--scaling and normalizing  $\tau$ according to the  braid equivalence in the given braid category. 
 In Table~1 we list some knot algebras with related knot invariants.

\begin{table}[h!]
  \begin{center}
\begin{tabular}{r l }
Knot algebra & Invariant \\
  \hline 
  Temperley--Lieb algebra & Jones polynomial \& bracket polynomial \\
  Iwahori--Hecke algebra & HOMFLYPT polynomial \\ 
	  BMW algebra & Kauffman polynomial \\ 
	  $B$--type \& affine Hecke algebras & Lambropoulou invariants\\
   Singular Hecke algebra &  Kauffman--Vogel \& Paris--Rabenda invariants \\
	Rook algebra & Alexander polynomial 
\end{tabular}
 \end{center}
  \caption{Examples of knot algebras.}
\end{table}

The inspiring example of framization is the so--called Yokonuma--Hecke algebra. Indeed, this algebra corresponds to a framization of the Iwahori--Hecke algebra. The Yokonuma--Hecke algebra appears in the field of group theory and was defined by Yokonuma \cite{yo} as a generalization of the Iwahori--Hecke algebra. More precisely, he considered the centralizer algebra of the permutation representation associated to 
any finite Chevalley group with respect to one maximal unipotent subgroup. Then, in analogy to the classical presentation of the  Iwahori--Hecke algebra, he found a presentation of the Yokonuma--Hecke algebra  by generators and relations \cite[Theorem 1]{yo}. 

\smallbreak
In \cite{ju4, juka} another presentation of the Yokonuma--Hecke algebra was established (Definition~\ref{yhdef}). By using this new presentation, the Yokonuma--Hecke algebra ${\rm Y}_{d, n}(u)$ was defined for any positive integers $n$ and $d$ and a fixed non--zero complex number $u$. 
 Then, ${\rm Y}_{d, n}(u)$ could be naturally viewed as a quotient of the framed braid group ${\mathcal F}_n$ or as an algebra obtained from the Iwahori--Hecke algebra ${\rm H}_n(u)$ by adding framing generators and by replacing the Hecke algebra quadratic  relation  by a  quadratic relation  which also involves intrinsically the framing generators (Eq.~\ref{quadr}). For $d=1$ the algebra ${\rm Y}_{d, n}(u)$ coincides with  the algebra ${\rm H}_n(u)$.

\smallbreak
Further, a Markov trace ${\rm tr}$ was constructed by the first author on the algebras ${\rm Y}_{d, n}(u)$, the `Juyumaya trace',  \cite{ju} with parameters $z, x_1, \ldots, x_{d-1}$, using an appropriate inductive linear basis. Parameter $z$ takes care of the braiding, while parameters $x_1, \ldots, x_{d-1}$ take care of the framing. For $d=1$ the trace  ${\rm tr}$ coincides with the well--known Ocneanu trace $\tau$ on the algebras ${\rm H}_n(u)$, from which the 2-variable Jones or HOMFLYPT polynomial (denoted here Homflypt)  for oriented links was extracted \cite{jo}.  

\smallbreak
Then, in order to obtain link invariants via the trace ${\rm tr}$, the `E--condition'  needed to be imposed on the framing parameters $x_1, \ldots, x_{d-1}$ (see Eq.~\ref{Econdition}) for re--scaling ${\rm tr}$ in order that it conforms with negative stabilization (cf.~\cite{jula5}). The trace  ${\rm tr}$  is the only known Markov trace on knot algebras that does not re--scale directly.  As it was shown by G\'erardin (cf. Appendix~\cite{jula5}) solutions of the `E--system' (\ref{Esystem}) are parametrized by the non--empty subsets of ${\mathbb Z}/d{\mathbb Z}$. 

\smallbreak
Subsequently, for any solution of the E--system, ${\rm tr}$ yielded an invariant for framed knots \cite{jula5}, for classical knots \cite{jula4} and  for singular knots \cite{jula3}. Moreover, the algebras ${\rm Y}_{d, n}(u)$ with the trace ${\rm tr}$ seem to relate naturally to the domain of transverse knots \cite{chjajukala}.
  All these invariants are still under investigation \cite{chla, chjajukala}, especially as to how they compare with the Homflypt polynomial.  In \cite{chla} it is shown that only in trivial cases the classical link invariants coincide with the Homflypt polynomial. By construction these invariants are at least as strong as the Homflypt polynomial, however computational evidence indicates that they are topologically equivalent.  This is not easy to see either by algebraic or by diagrammatic methods. In  \cite{chjajukala} some conjectures are stated in this direction. Still, we believe that it is  remarkable that  one can obtain invariants for all these different knot categories from a single algebra. 
	
\smallbreak	
The Yokonuma--Hecke algebras, equipped with a Markov trace, are interesting on their own right. Their representation theory has been studied thoroughly in \cite{th, chpa}. In particular, in \cite{chpa} a completely combinatorial approach is taken to the subject.
	
\smallbreak
All these results are presented in Sections~\ref{yhecke}, \ref{e} and \ref{invariants}.
 The above  comprise our motivation for constructing framizations of other knot algebras. 
 In this paper we  present possible framizations of most of the algebras listed in Table~1. Up to Section~\ref{trans} the paper is mostly a survey of results but from there on it continues as an announcement of new results, some by the authors, some  by the authors with co--authors and some by other authors. 

\smallbreak
In Section~\ref{tl} we  present three quotients of the Yokonuma--Hecke algebra ${\rm Y}_{d, n}(u)$ as possible framizations of the Temperley--Lieb algebra. For each one we present necessary and sufficient conditions for the trace   ${\rm tr}$ to pass through to the quotient algebra and we discuss related unoriented knot invariants \cite{gojukola, gojukolaf}. For the first quotient, the Yokonuma--Temperley--Lieb ${\rm YTL}_{d, n}(u)$, the ideal is similar to the one in the classical case. Then, as it turns out, the trace  ${\rm tr}$ passes to ${\rm YTL}_{d, n}(u)$ only if the  trace parameters $x_i$ are $d^{th}$  roots of unity. In this case we recover the Jones polynomial. See   \cite{gojukola}. The second  candidate, ${\rm FTL}_{d,n}(u)$, which we select as the `Framization of the Temperley--Lieb algebra', has the property that the conditions on the $x_i$'s so that   ${\rm tr}$ passes through to the quotient, include explicitely all solutions of the E--system mentioned above. Finally, the conditions on the trace parameters for the third candidate, the Complex Reflection Temperley--Lieb algebra ${\rm CTL}_{d,n}(u)$, involve only parameter $z$ and not the framing parameters $x_1, \ldots, x_{d-1}$. So, in order to obtain knot invariants from the algebra  ${\rm CTL}_{d,n}(u)$ we have to impose the E--condition (Eq.~\ref{Econdition}). It follows that the knot invariants we obtain coincide with those from ${\rm FTL}_{d,n}(u)$.  The main disadvantage of ${\rm CTL}_{d,n}(u)$ is that is is unnecessarily large for our purposes.

\smallbreak
In Section~\ref{b}  we  propose  framizations for the Hecke algebra of $B$--type, for the cyclotomic Hecke algebras of $B$--type, and for the generalized Hecke algebra of $B$--type, which is isomorphic to the affine Hecke algebra of $A$--type. These definitions were first given in \cite{jula6}. All these algebras are related to the knot theory of the solid torus \cite{la1,gela,la2}. These  $B$--type framizations are further studied in~\cite{chpa2}, where Markov traces are also constructed and a corresponding E--condition is given.

\smallbreak
In Section~\ref{bmw}  we  propose a framization of the Birman--Murakami--Wenzl or simply  BMW algebra  \cite{biwe,mu},  which is related to the Kauffman polynomial invariant of knots \cite{ka}.  This framization was introduced and further studied in \cite{jula6}. 

\smallbreak
Finally, in Section~\ref{sin} we propose a framization of the singular Hecke algebra \cite{para}, which is related to the invariants of Kauffman--Vogel  \cite{kavo} and Paris--Rabenda \cite{para}.  This was introduced in \cite{jula6}. 

\smallbreak
All these framization knot algebras are related to the framed braid group and they are of interest to algebraists. There are many more other knot algebras, such as other quotients of the classical braid group, quotients of the virtual braid group \cite{ka2,kala}, or the Rook algebra \cite{bry}  which is related to the Alexander polynomial. For all these one could construct appropriate framization counterparts.

\smallbreak
\noindent {\bf Acknowledgments} \ We would like to thank the Referee for the very careful reading and for the very interesting comments.

\section{Notations and Background}
 
\subsection{} 

Along the paper  the term algebra means a ${\mathbb C}$--associative  algebra with unity denoted $1$, where as usual 
${\mathbb C}$ denotes the field of complex numbers. Notice that ${\mathbb C}$ can be regarded as   included in the algebra as a central subalgebra. 
 We also denote by ${\mathbb C} G$ the group algebra of a group $G$. 

 \subsection{} 

The letters   $n$ and $d$ denote two positive integers. 
We denote by $S_n$ the symmetric group on $n$--symbols and by $s_i$ the elementary transposition $(i, \, i + 1)$. 
Let  $B_n$ be the {\it classical Artin braid group}. $B_n$ is presented  by the braiding generators $\sigma_1 , \ldots ,\sigma_{n-1}$ and the {\it braid relations}:
 \begin{equation}\label{braidrels}
 \begin{array}{rcll}
 \sigma_i\sigma_j & = & \sigma_j \sigma_i & \text{ for} \quad \vert i-j\vert > 1 \\  
\sigma_i\sigma_{j}\sigma_i  & = &    \sigma_{j}\sigma_i\sigma_{j} & \text{ for} \quad \vert i-j\vert =1.
\end{array}
\end{equation}
 The {\it framed braid group}  ${\mathcal F}_n$ is the group defined by adding to the above   presentation of $B_n$ the framing generators $t_1, \ldots , t_n$ and the following  relations:
 \begin{equation}\label{framerels}
 \begin{array}{rcll}
 t_it_j & = & t_j t_i & \text{ for} \quad 1\leq i,j \leq n \\  
t_j\sigma_i  & = &    \sigma_it_{s_i(j)} & \text{ for} \quad 1\leq i \leq n-1 \ \ \& \ \ 1\leq j \leq n,
\end{array}
\end{equation}
where $s_i(j)$ is the result of applying $s_i$ to $j$.  The {\it $d$--modular framed braid group}, denoted ${\mathcal F}_{d,n}$, is defined by adding to the above presentation of ${\mathcal F}_n$ the relations:
 \begin{equation}\label{modular}
 \begin{array}{rcll}
t_i^d & = & 1 & \text{ for}  \quad  1\leq i \leq n.
\end{array}
\end{equation}
We denote $C$ the infinite cyclic group and by $C_d$ the cyclic group of order $d$.  We have $C\cong {\mathbb Z}$ and $C_d \cong {\mathbb Z}/d{\mathbb Z}$. Further, if $t$ is a  generator of $C$, the group $C_d$ can be presented as  
$C_d = \left\langle t \,;\, t^d = 1\right\rangle$. From the above we have: ${\mathcal F}_n = C^n\rtimes B_n $ and ${\mathcal F}_{d,n} = C_d^n\rtimes B_n $. Finally, we shall denote the group: 
 \begin{equation}\label{cdn}
C_{d,n} := C_d^n\rtimes S_n.
\end{equation}

\subsection{} 

From now on we fix a non--zero complex number $u$ and a positive integer $d$. 

\begin{note}\rm 
One of the authors  does not agree with the denomination \lq Juyumaya trace\rq and the other author does not agree with the denomination \lq Lambropoulou invariants\rq. 
\end{note}

\section{The framization of the Hecke algebra of type $A$}\label{yhecke}

\subsection{}\label{heckeal} 

The {\it Iwahori--Hecke algebra of type $A$},  ${\rm H}_n(u)$, is  the algebra with  ${\mathbb C}$--linear basis $\{ h_w \, | \, w\in S_n \}$ and the following rules of multiplication: 
$$
h_{s_i}h_w =\left\{\begin{array}{ll}
h_{s_iw} & \text{for }  \  \ell(s_iw)>\ell(w) \\
 u h_{s_iw} + (u-1)h_w & \text{for } \ \ell(s_iw)<\ell(w) 
\end{array}\right.
$$
where $\ell$ is the usual length function on the symmetric group. Set $h_i := h_{s_i} $. 
 As usual we consider the presentation of the algebra ${\rm H}_n(u)$ by `braiding' generators $h_1, \ldots, h_{n-1}$, subject to the braid relations (\ref{braidrels}) together with the quadratic  relations:
\begin{equation}\label{hecke}
h_i^2 = (u-1)h_i + u.
\end{equation}
Note that ${\rm H}_n(1)$ coincides with the group algebra ${\mathbb C} S_n$.   
The above presentation of ${\rm H}_n(u)$ says that the mapping $\pi: \sigma_i \mapsto h_i$  defines an epimorphism from  ${\mathbb C} B_n$ onto ${\rm H}_n(u)$ whose kernel is generated 
by the  expressions: $\sigma^2_i  - (u-1)\sigma_i - u$.

\subsection{}\label{homflypt}

The Hecke algebra ${\rm H}_n(u)$ is the knot algebra used to define the 2--variable Jones or Homflypt polynomial  according to the Jones construction \cite{jo}. Namely, one uses the Markov braid equivalence on  $\cup_n B_n$,  comprising conjugation in the groups $B_n$ and positive and negative stabilization and destabilization ($\alpha {\sigma}_n \sim \alpha \sim \alpha \sigma_n^{-1}$; $\alpha \in B_n$), the map $\pi$ from $B_n$ to ${\rm H}_n(u)$ and the Ocneanu trace  defined on ${\rm H}_n(u)$: 

\begin{thm}[Theorem 5.1 \cite{jo}]\label{ocneanu}
Let $\zeta$ be an indeterminate over ${\mathbb C}$.  Then, there exists a unique ${\mathbb C}$--linear map  $\tau$ from the inductive limit of the family  $\{{\rm H}_n(u)\}_n$ to ${\mathbb C}(\zeta)$,  such that $\tau(1)=1$ and  satisfying the following rules  for all $a, b\in {\rm H}_n(u)$:
$$
\begin{array}{rcll}
\tau(ab) & = &\tau(ba)  \qquad &  \\
\tau(a h_n) & = & \zeta \, \tau(a) \qquad &  (\text{Markov  property} ).
\end{array}
$$
\end{thm}
 Diagrammatically, the Markov property of $\tau$  is illustrated on the left--hand side of Figure~\ref{tracerules}. 

\smallbreak
Then  $\tau$ has to be {\it normalized}, so that the closed braids $\widehat{\alpha}$ and $\widehat{\alpha {\sigma}_n}$ $(\alpha\in B_n)$ be assigned the same value of the invariant, and also {\it re--scaled}, so that the closed braids $\widehat{\alpha \sigma_n^{-1}}$ and $\widehat{\alpha {\sigma}_n}$ get also assigned the same value  of the invariant. 
 So, defining: 
$$
\lambda := \frac{1-u+ \zeta}{u \zeta} \qquad \text{and}  \qquad  C :=  \frac{1}{\zeta \sqrt{\lambda}}
$$
then the {\it Homflypt polynomial} $P = P(u, \lambda)$ of an oriented link $L$ is defined as follows
$$
P (\widehat{\alpha}) = C^{n-1}(\tau\circ \pi)(\alpha)
$$
where  $\alpha$ is a braid on $n$ strands whose closure $\widehat{\alpha}$ is isotopic to $L$.

\subsection{} 

In group theory the Hecke algebra has a natural generalization, the Yokonuma--Hecke algebra ${\rm Y}_{d,n}(u)$, see \cite{yo}.  In \cite{ju4, juka}  a new presentation was found  for this algebra. Subsequently, in \cite{ju} this new presentation was considered for the Yokonuma--Hecke algebra  and was proved that it supports a Markov trace.
We define ${\rm Y}_{d,n}(u)$  as follows.

\begin{defn}\label{yhdef} \rm 
The {\it Yokonuma--Hecke algebra of  type $A$}, ${\rm Y}_{d,n}(u)$, is the algebra defined by the braiding generators $g_1, \ldots, g_{n-1}$, the  framing generators $t_1, \ldots, t_n$ and the relations:
\begin{equation} \label{yh}
\begin{array}{rcll}
 g_ig_j & = & g_j g_i &  \text{ for} \quad \vert i-j\vert >1 \\  
g_ig_{j}g_i  & = & g_jg_ig_j &  \text{ for} \quad \vert i-j\vert =1\\
t_it_j & = & t_j t_i &  \text{ for} \quad 1\leq i,j \leq n \\  
t_jg_i  & = &    g_it_{s_i(j)} &  \text{ for} \quad1\leq i \leq n-1 \ \ \& \ \ 1\leq j \leq n \\
t_i^d & = & 1 &  \text{ for} \quad 1\leq i \leq n 
\end{array}
\end{equation}
where $s_i(j)$ denotes the result of applying $s_i$ to $j$, together with the 
following quadratic relations:
\begin{equation}\label{quadr}
g_i^2 = 1 + (u-1)e_i + (u-1)e_ig_i \qquad \text{ for all} \quad 1\leq i \leq n-1 
\end{equation}
where 
\begin{equation}\label{ei}
e_i := \frac{1}{d}\sum_{s=0}^{d-1}t_i^s t_{i+1}^{d-s}.
\end{equation}
\end{defn}

The elements $e_i$ are easily seen to be idempotents. Using this fact it follows from  the quadratic relations (\ref{quadr}) that the generators $g_i$ are invertible:
\begin{equation}\label{invrs}
g_i^{-1} = g_i + (u^{-1} - 1)\, e_i + (u^{-1} - 1)\, e_i \, g_i.
\end{equation}

From the above presentation it is clear that the algebra ${\rm Y}_{d,n}(u)$ is a quotient of the modular framed braid group algebra ${\mathbb C} {\mathcal F}_{d,n}$ under the quadratic relations (\ref{quadr}). This observation leads to diagrammatic interpretations for the elements in ${\rm Y}_{d,n}(u)$. For example,  the elements  $e_i$ (which are in ${\mathbb C} {\mathcal F}_{d,n}$ and in ${\rm Y}_{d,n}(u)$) can be represented as in Figure~\ref{ed1}.

\smallbreak
\begin{figure}[H]
\begin{picture}(320,60)

\put(0,38){$e_1 =$}
\put(37, 41){$1$}
\qbezier(36,40)(41,40)(46,40)
\put(36,30){$d$}
\qbezier(55,20)(45,40)(55,60) 
\qbezier(65,20)(65,40)(65,60)
\qbezier(80,20)(80,40)(80,60)
\qbezier(95,20)(95,40)(95,60)
\put(105,40){$+$}
\qbezier(125,20)(125,40)(125,60)
\qbezier(140,20)(140,40)(140,60)
\qbezier(155,20)(155,40)(155,60)
\put(165,40){$+$}
\qbezier(185,20)(185,40)(185,60)
\qbezier(200,20)(200,40)(200,60)
\qbezier(215,20)(215,40)(215,60)
\put(225,40){$+$}
\put(240,40){$\cdots$}
\put(260,40){$+$}

\qbezier(280,20)(280,40)(280,60)
\qbezier(295,20)(295,40)(295,60)
\qbezier(310,20)(310,40)(310,60)
\qbezier(320,20)(330,48)(320,60) 
\put(62,65){\tiny{$0$}}
\put(77,65){\tiny{$0$}}
\put(92,65){\tiny{$0$}}
\put(122,65){\tiny{$1$}}
\put(130,65){\tiny{$d-1$}}
\put(153,65){\tiny{$0$}}
\put(182,65){\tiny{$2$}}
\put(190,65){\tiny{$d-2$}}
\put(213,65){\tiny{$0$}}
\put(270,65){\tiny{$d-1$}}
\put(292,65){\tiny{$1$}}
\put(307,65){\tiny{$0$}}

\end{picture}
\caption{The element $e_1 \in {\mathbb C} {\mathcal F}_{d,3}$.} 
\label{ed1}
\end{figure}

\noindent Also, Eq.~\ref{invrs}, which is equivalent to the quadratic relation, is illustrated in Figure~\ref{g1invrs}.

\smallbreak
\begin{figure}[H]
\setlength{\unitlength}{.8pt}
\begin{picture}(400,150)

\qbezier(-50,100)(-51,105)(-50,110)
\qbezier(-30,100)(-29,105)(-30,110)

\qbezier(-50,110)(-50,114)(-45,117)
\qbezier(-35,123)(-30,126)(-30,130)
\qbezier(-30,110)(-30,115)(-40,120)
\qbezier(-40,120)(-50,125)(-50,130)

\qbezier(-50,130)(-49,135)(-50,140)
\qbezier(-30,130)(-29,135)(-30,140)
\qbezier(-10,100)(-10,120)(-10,140)
\put(0,117){$=$}

\qbezier(20,100)(18,105)(20,110)
\qbezier(40,100)(41,105)(40,110)
\qbezier(20,110)(20,115)(30,120)
\qbezier(30,120)(40,125)(40,130)
\qbezier(40,110)(40,114)(35,117)
\qbezier(25,123)(20,126)(20,130)
\qbezier(20,130)(19.5,135)(20,140)
\qbezier(40,130)(41,135)(40,140)

\qbezier(57,100)(57,120)(57,140)

\put(64,117){\small{$+\frac{u^{-1}-1}{d}$}}
\qbezier(117,100)(107,120)(117,140) 
\qbezier(130,100)(130,120)(130,140)
\qbezier(148,100)(148,120)(148,140)
\qbezier(166,100)(166,120)(166,140)
\put(180,117){$+$}
\qbezier(203,100)(203,120)(203,140)
\qbezier(221,100)(221,120)(221,140)
\qbezier(239,100)(239,120)(239,140)
\put(250,117){$+$}
\qbezier(273,100)(273,120)(273,140)
\qbezier(291,100)(291,120)(291,140)
\qbezier(309,100)(309,120)(309,140)
\put(323,117){$+\cdots +$}
\qbezier(380,100)(380,120)(380,140)
\qbezier(398,100)(398,120)(398,140)
\qbezier(416,100)(416,120)(416,140)
\qbezier(425,100)(435,120)(425,140) 

\put(-3,37){\small{$+\frac{u^{-1}-1}{d}$}}

\qbezier(50,20)(40,40)(50,60) 

\qbezier(60,50)(59,55)(60,60)
\qbezier(80,50)(81,55)(80,60)

\qbezier(60,20)(59,25)(60,30)
\qbezier(80,20)(81,25)(80,30)

\qbezier(60,30)(60,35)(70,40)
\qbezier(70,40)(80,45)(80,50)
\qbezier(80,30)(80,34)(75,37)
\qbezier(65,43)(60,46)(60,50)

\qbezier(95,20)(95,40)(95,60)

\put(115,37){$+$}

\qbezier(150,50)(149,55)(150,60)
\qbezier(170,50)(171,55)(170,60)

\qbezier(150,20)(149,25)(150,30)
\qbezier(170,20)(171,25)(170,30)

\qbezier(150,30)(150,35)(160,40)
\qbezier(160,40)(170,45)(170,50)
\qbezier(170,30)(170,34)(165,37)
\qbezier(155,43)(150,46)(150,50)
\qbezier(185,20)(185,40)(185,60)

\put(205,37){$+$}

\qbezier(230,50)(229,55)(230,60)
\qbezier(250,50)(251,55)(250,60)

\qbezier(230,20)(229,25)(230,30)
\qbezier(250,20)(251,25)(250,30)

\qbezier(230,30)(230,35)(240,40)
\qbezier(240,40)(250,45)(250,50)
\qbezier(250,30)(250,34)(245,37)
\qbezier(235,43)(230,46)(230,50)
\qbezier(265,20)(265,40)(265,60)
\put(285,37){$+\cdots +$}

\qbezier(340,50)(339,55)(340,60)
\qbezier(360,50)(361,55)(360,60)

\qbezier(340,20)(339,25)(340,30)
\qbezier(360,20)(361,25)(360,30)

\qbezier(340,30)(340,35)(350,40)
\qbezier(350,40)(360,45)(360,50)
\qbezier(360,30)(360,34)(355,37)
\qbezier(345,43)(340,46)(340,50)

\qbezier(375,20)(375,40)(375,60)

\qbezier(385,20)(395,40)(385,60) 
\put(-50,142){\tiny{$0$}}
\put(-32,142){\tiny{$0$}}
\put(-11,142){\tiny{$0$}}

\put(18,142){{\tiny $0$}}
\put(38,142){\tiny{$0$}}
\put(55,142){\tiny{$0$}}

\put(127,142){\tiny{$0$}}
\put(146,142){\tiny{$0$}}
\put(164,142){\tiny{$0$}}

\put(200,142){\tiny{$1$}}
\put(213,142){\tiny{$d\!-\!1$}}
\put(237,142){\tiny{$0$}}

\put(270,142){\tiny{$2$}}
\put(282,142){\tiny{$d\!-\!2$}}
\put(307,142){\tiny{$0$}}

\put(366,142){\tiny{$d\!-\!1$}}
\put(395,142){\tiny{$1$}}
\put(414,142){\tiny{$0$}}

\put(58,63){\tiny{$0$}}
\put(77,63){\tiny{$0$}}
\put(93,63){\tiny{$0$}}

\put(148,63){\tiny{$1$}}
\put(160,63){\tiny{$d\!-\!1$}}
\put(183,63){\tiny{$0$}}

\put(228,63){\tiny{$2$}}
\put(240,63){\tiny{$d\!-\!2$}}
\put(264,63){\tiny{$0$}}

\put(331,63){\tiny{$d\!-\!1$}}
\put(358,63){\tiny{$1$}}
\put(372,63){\tiny{$0$}}
\end{picture}
\caption{The element $g_1^{-1} \in {\rm Y}_{d,3}(u)$.}
\label{g1invrs}
\end{figure}

\subsection{}

A consequence of the  above definition is that every word in the defining generators of  ${\rm Y}_{d,n}(u)$ can be written in the {\it split form} $t_1^{a_1}\ldots t_n^{a_n} \, g$, where the $a_i$'s are integers modulo $d$ and $g$ is a word in the $g_i$'s. Since the $g_i$'s satisfy the braid relations we have that, if $w=s_{i_1}\ldots s_{i_m} \in S_n$ is a reduced expression, then the following element $g_w:= g_{i_1}\ldots g_{i_{m}}$ is well--defined. 
 In \cite{ju4} it is proved that the multiplication rules in ${\rm Y}_{d,n}(u)$ are governed by the 
 group $C_{d,n}$. In fact, the multiplication rules among the framing generators and between the  framing generators  and  the braiding generators are the same multiplication rules  as in the group $C_{d,n}$.  For the multiplication among the braiding generators we have:
 $$
g_ig_w =\left\{\begin{array}{ll}
g_{s_iw} & \text{for} \ \ell (s_iw)>\ell (w) \\
 g_{s_iw} + (u-1)e_ig_{s_iw}  +  (u-1)e_ig_w & \text{for} \ \ell (s_iw)< \ell (w) 
\end{array}\right.
$$

Notice now that ${\rm Y}_{d,n}(1) = {\mathbb C}C_{d,n}$. This says that ${\rm Y}_{d,n}(u)$ is essentially obtained from  ${\rm H}_n(u)$ by  adding framing generators, since  ${\rm H}_n(1) = {\mathbb C}S_n$. For this reason we shall call the Yokonuma--Hecke algebra ${\rm Y}_{d,n}(u)$ a  {\it framization of the Iwahori--Hecke algebra} ${\rm H}_n(u)$.  The key point  in this framization is the quadratic relation (\ref{quadr}), which is considered as the framization of the Hecke algebra quadratic relation. 

\smallbreak
The representation theory  of the Yokonuma--Hecke algebra has been  studied in \cite{th} and \cite{chpa}. Finally, it is worth mentioning that the notions of the modular framed braid group ${\mathcal F}_{d,n}$ as well as of the Yokonuma--Hecke algebra  ${\rm Y}_{d,n}(u)$ have been extended to constructions  on the $p$--adic level  \cite{jula1, jula2} and the adelic level \cite{jula4}.

\subsection{}

\begin{figure}[H]
\begin{picture}(330,70)
\put(-23,36){${\rm tr}$}
\qbezier(0,0)(-9,37)(0,74)  
\put(13,36){$a$}
\qbezier(8,0)(8,12)(8,24)
\qbezier(8,54)(8,62)(8,70)
\qbezier(22,54)(22,62)(22,70)
\qbezier(50,24)(52,47)(50,70)
\qbezier(30,0)(29,2)(30,4)
\qbezier(50,0)(51,2)(50,4)
\qbezier(58,0)(67,37)(58,74)  
\put(68,36){$=$}
\put(78,36){$z\,{\rm tr}$}
\qbezier(107,0)(96,37)(107,74)  
\qbezier(0,24)(0,39)(0,54)
\qbezier(30,24)(30,39)(30,54)

\qbezier(0,54)(15,54)(30,54)
\qbezier(0,24)(15,24)(30,24)
\qbezier(30,4)(30,9)(40,14)
\qbezier(40,14)(50,19)(50,24)
\qbezier(50,4)(50,8)(45,11)
\qbezier(35,17)(30,20)(30,24)
\qbezier(110,24)(110,39)(110,54)
\qbezier(140,24)(140,39)(140,54)

\qbezier(110,54)(125,54)(140,54)
\qbezier(110,24)(125,24)(140,24)
\qbezier(118,0)(118,12)(118,24)
\qbezier(118,54)(118,62)(118,70)
\qbezier(132,54)(132,62)(132,70)
\qbezier(132,0)(132,12)(132,24)
\put(122,36){$a$}
\qbezier(146,0)(156,37)(146,74) 
\put(160,36){$,$}
\put(175,36){${\rm tr}$}
\qbezier(197,0)(188,37)(197,74)  
\qbezier(202,24)(202,39)(202,54)
\qbezier(232,24)(231,39)(232,54)

\qbezier(202,54)(217,54)(232,54)
\qbezier(202,24)(217,24)(232,24)
\qbezier(210,0)(210,12)(210,24)
\qbezier(224,54)(224,62)(224,70)
\qbezier(210,54)(210,62)(210,70)
\qbezier(224,0)(224,12)(224,24)

\qbezier(240,0)(240,34)(240,68)
\put(238, 70){\tiny{ $m$}}
\put(215,36){$a$}
\qbezier(253,0)(261,37)(253,74) 
\put(261,36){$=$}
\put(272,36){$x_m {\rm tr}$}
\qbezier(306,0)(298,37)(306,74)  
\qbezier(310,24)(310,39)(310,54)
\qbezier(340,24)(340,39)(340,54)

\qbezier(310,54)(325,54)(340,54)
\qbezier(310,24)(325,24)(340,24)
\qbezier(318,0)(318,12)(318,24)
\qbezier(318,54)(318,62)(318,70)
\qbezier(332,54)(332,62)(332,70)
\qbezier(332,0)(332,12)(332,24)

\put(321,36){$a$}
\qbezier(345,0)(353,37)(345,74) 
\end{picture}
\caption{Topological interpretations of the trace rules}
\label{tracerules}
\end{figure}
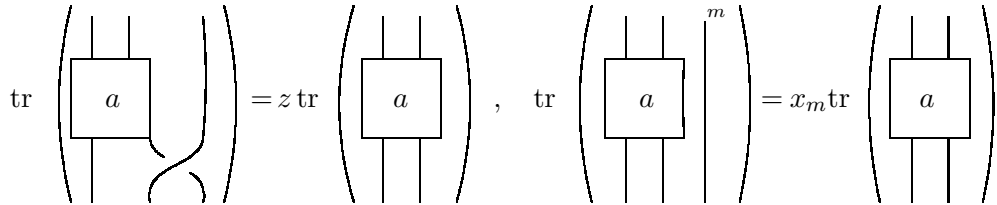

Another crucial property of the Yokonuma--Hecke algebra is that is supports a Markov trace:

\begin{thm}[Theorem 12 \cite{ju}]\label{trace}
 Let $z$, $x_1, \ldots ,x_{d-1}$ be indeterminates over ${\mathbb C}$. Then, there exists a unique ${\mathbb C}$--linear map 
 ${\rm tr}$ from the inductive limit of the family  $\{{\rm Y}_{d,n}(u)\}_n$ to ${\mathbb C}(z, x_1, \ldots ,x_{d-1})$, 
 such that ${\rm tr}(1)=1$ and  satisfying the following rules  for all $a, b\in {\rm Y}_{d,n}(u)$:
$$
\begin{array}{crcll}
(1) & {\rm tr}(ab) & = & {\rm tr}(ba)  \qquad &  \\
(2) & {\rm tr}(a \, g_n) & = & z \, {\rm tr}(a) \qquad &  (\text{Markov  property} ) \\
(3) & {\rm tr}(a \, t_{n+1}^m) & = & x_m {\rm tr}(a) \qquad & (1\leq m\leq d-1).
\end{array}
$$
\end{thm}
The topological interpretations for rules (2) and (3) are given in Figure~\ref{tracerules}.

The trace ${\rm tr}$ shall be called the {\it Juyumaya trace}. Note that for $d=1$, the algebra ${\rm Y}_{1,n}(u)$ is the Iwahori--Hecke algebra ${\rm H}_n(u)$ and the Juyumaya  trace coincides with the Ocneanu trace with parameter $z$. Note also that the trace  ${\rm tr}$ lifts to the $p$--adic level \cite{jula1, jula5} and to the adelic level  \cite{jula4}.

\subsection{}

A Markov trace is a key ingredient for constructing knot invariants from an algebra. The first key requirement is that there is a representation of some braid category to the algebra.   Another key requirement is to have a Markov--type braid equivalence for the braid category, corresponding to isotopy in the related knot category. The last key ingredient is to re--scale and normalize the trace according to the given braid equivalence. 

\smallbreak
 In the case of the Yokonuma--Hecke algebras we have a natural connection with the framed braid category and with the classical braid category. More precisely, the defining relations of ${\rm Y}_{d,n}(u)$ yield two natural representations. One  of the framed braid group:
\begin{equation}\label{gamma}
\begin{array}{rccc}
 \gamma : {\mathbb C} {\mathcal  F}_{n} & \longrightarrow &  {\rm Y}_{d,n}(u) & \\
\sigma_i  &  \mapsto & g_i &  \\
t_j^s   & \mapsto & t_j^{s( {\rm mod}\, d)}
\end{array}
\end{equation}
 and another  of  the classical braid group: 
\begin{equation}\label{delta}
\begin{array}{rccc}
 \delta : {\mathbb C} B_{n}  & \longrightarrow &  {\rm Y}_{d,n}(u) & \\
\sigma_i  &  \mapsto & g_i &  
\end{array}
\end{equation}
The map $\delta$ can be viewed as the composition of the map $\gamma$ with the natural injection of $B_n$ into ${\mathcal  F}_{n}$, whereby  $\sigma_i \mapsto \sigma_i$, and classical braids are considered to have all framings zero. 

\smallbreak
In the sequel we will sometimes identify a braid $\alpha$ with its image through $\gamma$ or $\delta$.

\section{The E--system}\label{e}

The framed braid equivalence comprises conjugation in the groups ${\mathcal F}_n$ and positive and negative stabilization and destabilization (see for example \cite{ks}). So, with the trace ${\rm tr}$ in hand, in \cite{jula5} the authors tried to obtain topological invariants for framed links after the method of V.F.R.~Jones \cite{jo}. 
 This meant that ${\rm tr}$ would have to be normalized, so that the closed framed braids $\widehat{\alpha}$ and  $\widehat{\alpha\, {\sigma}_n}$ $(\alpha\in {\mathcal F}_n)$ get assigned the same value of the invariant, and also re--scaled, so that the closed framed braids $\widehat{\alpha\, \sigma_n^{-1}}$ and $\widehat{\alpha\, {\sigma}_n}$  $(\alpha\in {\mathcal F}_n)$ get also assigned the same value  of the invariant.  However, as it turned out, ${\rm tr}(\alpha\, g_n^{-1})$ does not factor through ${\rm tr}(\alpha)$. That is, remarkably: 
\begin{equation}\label{nofactor}
{\rm tr}(\alpha\, g_n^{-1})  
 \neq {\rm tr}(g_n^{-1}) {\rm tr}(\alpha), 
\end{equation}
since, from (\ref{invrs}), we deduce that 
${\rm tr} (\alpha g_n^{-1})={\rm tr}(\alpha\, g_n) + (u^{-1} - 1)\, {\rm tr}(\alpha\, e_n) +(u^{-1} - 1)\, {\rm tr}(\alpha\, e_n g_n)$
  and, although ${\rm tr}(\alpha\, e_n g_n) = z\, {\rm tr}(\alpha)$, the term  
 ${\rm tr}(\alpha\, e_{n})$ does not factor through ${\rm tr}(\alpha)$.   Forcing the {\it {\rm E}--condition}: 
\begin{equation}\label{Econdition}
{\rm tr}(\alpha\, e_{n}) = {\rm tr}(e_{n})\, {\rm tr}(\alpha) \qquad (a \in {\rm Y}_{d,n}(u))
\end{equation}
implies equivalently that the trace parameters $x_1, \ldots, x_{d-1}$ have to satisfy the  E--{\it system} \cite{jula5},  the following non--linear system of equations in ${\mathbb C}$:
\begin{equation}\label{Esystem}
E^{(m)} = {x}_m E \qquad (1\leq m \leq d-1)
\end{equation}
where
$$
E = E^{(0)} := {\rm tr}(e_{i}) = \frac{1}{d}\sum_{s=0}^{d-1}{x}_{s}{x}_{d-s} 
 \qquad \mbox{and} \qquad 
E^{(m)} :=\frac{1}{d} \sum_{s=0}^{d-1}{x}_{m+s}{x}_{d-s},
$$
where the sub--indices on the ${x}_j$'s are regarded modulo $d$ and ${x}_0:=1$. 

As it was shown  by P.~G\'{e}rardin (\cite[Appendix]{jula5}), the solutions of the E--system are parametrized by the non--empty subsets of ${\mathbb Z}/d{\mathbb Z}$. 

\smallbreak
It is worth noting that the solutions of the E--system can be interpreted as a generalization of the Ramanujan sum. Indeed, by taking the subset $R$ of  $\mathbb{Z}/d{\mathbb{Z}}$ comprising the numbers coprime to $d$, then the solution parametrized by $R$ is, up to the factor $\vert R\vert$,  the Ramanujan sum $c_d(k)$ (see \cite{ra}). 

\smallbreak
It is also worth mentioning that solutions of the E--system lift to solutions on the $p$--adic level  \cite{jula1,jula5} and on the adelic level \cite{jula4}.

\section{Knot invariants from the Yokonuma--Hecke algebras}\label{invariants}

The Yokonuma--Hecke algebras have been used for constructing invariants for framed knots \cite{jula5}, for classical knots \cite{jula4} and for singular knots \cite{jula3}.  These invariants qualify the  algebra  ${\rm Y}_{d,n}(u)$ as a knot algebra and they comprise our main motivation for extending the notion of framization to other known knot algebras. We shall now recall briefly the construction of  these invariants.  

\subsection{}

Let $X_D= ( {\rm x}_1, \ldots , {\rm x}_{d-1} )$ be a solution of the E--system parametrized by the non--empty subset $D$ of ${\mathbb Z}/d{\mathbb Z}$. 

\begin{defn}[Definition~3 \cite{chla}] \rm
The trace map ${\rm tr}_D$ defined as the trace ${\rm tr}$  with the parameters $x_i$ specialized to the values ${\rm x}_i$, shall be called the \emph{specialized Juyumaya trace}  with parameter $z$.
\end{defn}
 Note that for $d=1$ the trace ${\rm tr}_D$  coincides with the trace ${\rm tr}$ (and with  the Ocneanu trace with parameter $z$). As it turns out \cite{jula4, jula5},  
$$
 E_D := {\rm tr}_D(e_{i})  = \frac{ 1}{\vert D\vert}
$$ 
where $\vert D\vert$ is the cardinal of the subset $D$.

\subsection{}\label{framedknots}

Let ${\mathcal L}_f$ denote the set of oriented  framed links. From the above, and re--scaling $g_i$ to $\sqrt{\lambda_D}\, g_i$, so that  ${\rm tr}_D(g_i^{-1}) =\lambda_D \, z$, where:
\begin{equation}\label{Dlambda}
\lambda_D := \frac{z +(1-u)E_D}{uz} = \frac{\vert D \vert z +1 - u}{\vert D \vert uz}
\end{equation}
 we have the following (mapping $\sigma_i\mapsto \sqrt{\lambda_D}\, g_i$):

\begin{thm}[\cite{jula5}]\label{framedinv} 
 Given a solution $X_D$ of the E--system, for any framed braid $\alpha \in {\mathcal  F}_{n}$ we define for the framed link $\widehat{\alpha} \in {\mathcal L}_f$:
$$
\Gamma_D (\widehat{\alpha}) = P_D^{n-1} (\sqrt{\lambda_D})^{\epsilon(\alpha)} \left({\rm tr}_D \circ \gamma \right)(\alpha)
$$
where 
$P_D = \frac{1}{z \sqrt{\lambda_D}}$, 
$\epsilon(\alpha)$ is the algebraic sum of the exponents of the $\sigma_i$'s in $\alpha$ and $\gamma$ the epimorphism (\ref{gamma}).
 Then the map  $\Gamma_D (u,\lambda_D)$  is a 2--variable isotopy invariant of oriented framed links.
\end{thm}

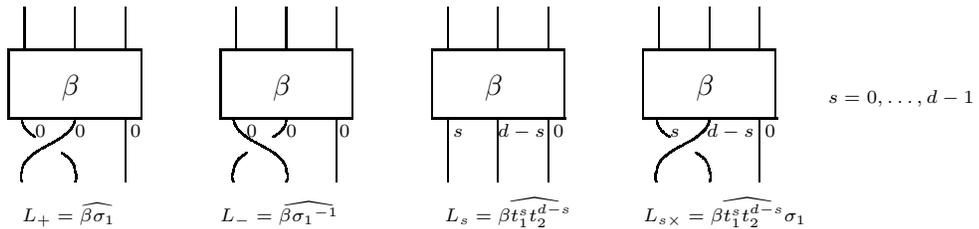
\begin{figure}[H]
 
\begin{picture}(330,80)

\put(20,53){$\beta$}

\qbezier(6,70)(6,78)(6,86) 
\qbezier(25,70)(25,78)(25,86)
\qbezier(44,70)(44,78)(44,86)

\qbezier(5,20)(4,22)(5,24) 
\qbezier(25,20)(26,22)(25,24) 

\qbezier(0,44)(0,59)(0,70)
\qbezier(50,44)(50,59)(50,70)

\qbezier(0,70)(25,70)(50,70)
\qbezier(0,44)(25,44)(50,44)
\qbezier(5,24)(5,29)(15,34)
\qbezier(15,34)(25,39)(25,44)
\qbezier(25,24)(25,28)(20,31)
\qbezier(10,37)(5,40)(5,44)
\qbezier(44,20)(44,32)(44,44)


\put(100,53){$\beta$}

\qbezier(86,70)(86,78)(86,86) 
\qbezier(105,70)(105,78)(105,86)
\qbezier(124,70)(124,78)(124,86)

\qbezier(85,20)(84,22)(85,24) 
\qbezier(105,20)(106,22)(105,24) 

\qbezier(80,44)(80,59)(80,70)
\qbezier(130,44)(130,59)(130,70)

\qbezier(80,70)(105,70)(130,70)
\qbezier(80,44)(105,44)(130,44)

\qbezier(85,24)(85,28)(90,31)
\qbezier(100,37)(105,40)(105,44)
\qbezier(105,24)(105,29)(95,34)
\qbezier(95,34)(85,39)(85,44)

\qbezier(124,20)(124,32)(124,44)


\put(180,53){$\beta$}

\qbezier(166,70)(166,78)(166,86) 
\qbezier(185,70)(185,78)(185,86)
\qbezier(204,70)(204,78)(204,86)
\qbezier(160,44)(160,59)(160,70)
\qbezier(210,44)(210,59)(210,70)

\qbezier(160,70)(185,70)(210,70)
\qbezier(160,44)(215,44)(210,44)

\qbezier(166,20)(166,32)(166,44) 
\qbezier(185,20)(185,32)(185,44)
\qbezier(204,20)(204,32)(204,44)

\put(260,53){$\beta$}

\qbezier(246,70)(246,78)(246,86) 
\qbezier(265,70)(265,78)(265,86)
\qbezier(284,70)(284,78)(284,86)

\qbezier(245,20)(244,22)(245,24) 
\qbezier(265,20)(266,22)(265,24) 

\qbezier(240,44)(240,59)(240,70)
\qbezier(290,44)(290,59)(290,70)

\qbezier(240,70)(265,70)(290,70)
\qbezier(240,44)(295,44)(290,44)
\qbezier(245,24)(245,29)(255,34)
\qbezier(255,34)(265,39)(265,44)
\qbezier(265,24)(265,28)(260,31)
\qbezier(250,37)(245,40)(245,44)
\qbezier(284,20)(284,32)(284,44)
\put(10,37){\tiny{$0$}}
\put(25,37){\tiny{$0$}}
\put(46,37){\tiny{$0$}}

\put(90,37){\tiny{$0$}}
\put(105,37){\tiny{$0$}}
\put(125,37){\tiny{$0$}}

\put(168,37){\tiny{$s$}}
\put(185,37){\tiny{$d-s$}}
\put(206,37){\tiny{$0$}}

\put(250,37){\tiny{$s$}}
\put(264,37){\tiny{$d-s$}}
\put(286,37){\tiny{$0$}}

\put(310,50){\tiny{$s=0, \ldots ,d-1$}}
\put(5,5){\tiny{$L_{+}=\widehat{\beta \sigma_1}$}}
\put(80,5){\tiny{$L_{-}=\widehat{\beta {\sigma_1}^{-1}}$}}
\put(165,5){\tiny{$L_s=\widehat{\beta t_1^st_2^{d-s}}$}}
\put(240,5){\tiny{$L_{s\times}=\widehat{\beta t_1^st_2^{d-s}}\sigma_1$}}

\end{picture}
\caption{The framed links $L_+$, $L_-$, $L_s$ and  $L_{s\times}$.}
\label{skeinffig}
\end{figure}

 Further, in \cite{jula5} a skein relation has been found for the invariant  $\Gamma_D (u,\lambda_D)$ involving the braiding and the framing generators:
\begin{equation}\label{skeinf}
 \sqrt{\lambda} \,\Gamma_D (L_{-})
= \frac{1}{\sqrt{\lambda}} \, \Gamma_D (L_{+})  + \frac{u^ {-1}-1}{d} \sum_{s=0}^{d-1}\Gamma_D (L_s)
+ \frac{u^ {-1}-1}{d\sqrt{\lambda}} \sum_{s=0}^{d-1}\Gamma_D (L_{s\times})
\end{equation}
where $L_{+}$, $L_{-}$, $L_s$ and $L_{s\times}$, $s=0,\ldots,d-1$, are diagrams of oriented framed links, which are all identical except near one crossing, where they differ by the ways indicated in Figure~\ref{skeinffig}. 

\smallbreak

 Finally,  a $p$--adic (resp. adelic) invariant for oriented framed links, $\Gamma_{p^{\infty}}$, has been defined through the invariants $\Gamma_D$. For more details see \cite{jula5, jula4}.

\begin{note} \rm It is logical that one should try to extract from the invariants for framed links invariants of 3--manifolds. But for this, one needs to start with invariants of {\it unoriented} framed links. This can be achieved from our constructions so far by taking appropriate quotients of the Yokonuma--Hecke algebras, analogous to the Temperley--Lieb algebras. This is presented in Section~\ref{tl}. 
\end{note}

\subsection{}\label{classical}

Let  ${\mathcal L}$ denote the set of oriented classical links. A classical link may be viewed as a framed link with all framings zero. By the mapping (\ref{delta}) of the classical braid group $B_n$ in ${\rm Y}_{d,n}(u)$, by the classical Markov braid equivalence and using the construction and notations above we obtain invariants for classical links, where the $t_j$'s are treated as formal generators \cite{jula4}. Namely:

\begin{thm}[\cite{jula4}]\label{classicinv} 
 Given a solution $X_D$ of the E--system, for any braid $\alpha \in B_{n}$ we define for the link $\widehat{\alpha} \in {\mathcal L}$:
$$
\Delta_D (\widehat{\alpha}) = P_D^{n-1} (\sqrt{\lambda_D})^{\epsilon(\alpha)} \left({\rm tr}_D \circ \delta \right)(\alpha)
$$
where $P_D, \lambda_D$ are as defined in \S\ref{framedknots}, $\epsilon(\alpha)$ is the algebraic sum of the exponents of the $\sigma_i$'s in $\alpha$ and $\delta$ is the natural algebra homomorphism (\ref{delta}).
 Then the map $\Delta_D (u,\lambda_D)$ is a $2$--variable isotopy invariant of oriented links. 
\end{thm}

The invariant $\Delta_D (u,\lambda_D)$ can be viewed as the restriction of $\Gamma_D (u,\lambda_D)$ on the set of framed links with all framings zero. 

\begin{figure}[H]
\begin{center}
 
\begin{picture}(230,80)
\put(29,80){\line(-1,-1){28}}
\put(17,63){\line(1,-1){11}}
\put(0,80){\line(1,-1){11}}

\put(28,52){\vector(-1,-1){28}}
\put(18,35){\vector(1,-1){11}}
\put(0,52){\line(1,-1){11}}


\put(120,80){\vector(-1,-1){50}}
\put(98,52){\vector(1,-1){22}}
\put(70,80){\line(1,-1){22}}


\put(160,80){\vector(0,-1){50}}
\put(190, 80){\vector(0,-1){50}}


\put(258,58){\line(1,1){22}}
\put(252,52){\vector(-1,-1){22}}
\put(230,80){\vector(1,-1){50}}

\put(12,10){$L_{+ +}$}
\put(90,10){$L_{+}$}
\put(170,10){$L_{0}$}
\put(250,10){$L_{-}$}

\end{picture}
\caption{The classical links $L_{++}$, $L_+$, $L_0$ and $L_-$.}\label{skeincfig}
\end{center}
\end{figure}
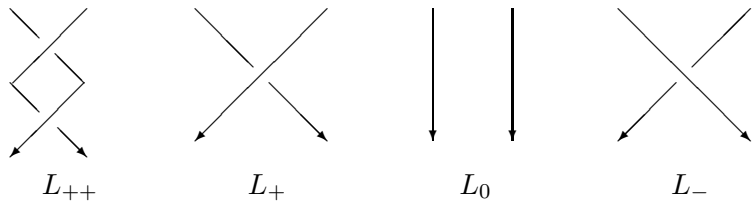

\smallbreak
The invariants $\Delta_D(u,\lambda_D)$ need to be compared with known invariants of classical links, especially with the Homflypt polynomial $P(u,\lambda)$, recall \S\ref{homflypt}. This is not easy to do on the algebraic level as there are no algebra homomorphisms connecting the algebras and the traces \cite{chla}. 
 Further, the skein relation of the invariant $\Delta_D(u,\lambda_D)$ has no topological interpretation in the case of classical links.  This makes the comparison very difficult also using diagrammatic methods. It is worth noting at this point  that in ${\rm Y}_{d,n}(u)$ a `closed' cubic relation is satisfied \cite{jula2, jula4}, closed in the sense of involving only the braiding generators, which is of minimal degree. Namely: 
\begin{equation}\label{cubic}
g_i^3 = u g_i^2 + g_i - u.
\end{equation}
 The cubic relation gives rise to a cubic skein relation for the invariant $\Delta_D(u,\lambda_D)$, involving only the braiding generators \cite{jula4}:
\begin{equation}\label{skeinc}
\sqrt{\lambda_D}\, \Delta_D (L_{-}) =
- \frac{1}{u \, \lambda_D}\, \Delta_D (L_{++})  + \frac{1}{\sqrt{\lambda_D}}\, \Delta_D (L_{+}) + \frac{1}{u}\, \Delta_D (L_{0})
\end{equation}
where the links $L_{++}$, $L_+$, $L_0$ and $L_-$ have identical diagrams, except for a region where they differ in the manners illustrated in Figure~\ref{skeincfig}. However, the above skein relation is not sufficient for determining the invariant $\Delta_D(u,z)$ diagrammatically with a simple set of initial conditions.

\smallbreak

 As it turns out, the cubic relation factors to the quadratic relation of the Iwahori--Hecke algebra ${\rm H}_n(u)$: 
\begin{equation}\label{cubicanalysis}
g_i^3 - u g_i^2 - g_i + u = (g_i - 1) \big(g_i^2 - (u - 1) g_i - u\big).
\end{equation}
Unfortunately, this factoring does not give information about the comparison of the invariants. 
 
\smallbreak
From our construction it follows that we have, at least,  recovered the Homflypt polynomial. 
 In \cite{chla} it is shown that for generic values of the parameters $u,z$ the invariants $\Delta_D(u,\lambda_D)$ do not coincide with the Homflypt polynomial except in the trivial cases $u=1$ or $E=1$. More precisely, for $E=1$ an algebra homomorphism $h : {\rm Y}_{d,n}(u) \longrightarrow {\rm H}_n(u)$ can be defined and the composition $\tau \circ h$ is a Markov trace on ${\rm Y}_{d,n}(u)$ which takes the same values as the specialized Juyumaya trace ${\rm tr}_D$ whereby the $x_i$'s are specialized to $d^{th}$ roots of unity \cite[Appendix]{jula5}. In this case we also obtain $\vert D \vert = 1$. For details see \cite[\S 3]{chla}. 

\smallbreak
Yet, computational data \cite{chjajukala} indicate that these invariants do not distinguish more knot pairs than the Homflypt polynomial, so they are probably topologically equivalent to the Homflypt polynomial. A complete answer to this question is still under investigation. In \cite{chjajukala} some conjectures are formulated and tackled in this direction. Also, it is proved that the specialized Juyumaya trace can be computed for classical braids with rule~(3) replaced by another rule involving the elements $e_i$.

\subsection{}

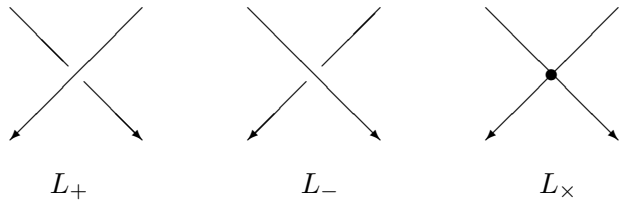
\begin{figure}
\begin{center}
 
\begin{picture}(230,80)
\put(50,80){\vector(-1,-1){50}}
\put(28,52){\vector(1,-1){22}}
\put(0,80){\line(1,-1){22}}

\put(118,58){\line(1,1){22}}
\put(112,52){\vector(-1,-1){22}}
\put(90,80){\vector(1,-1){50}}

\put(230,80){\vector(-1,-1){50}}
\put(180,80){\vector(1,-1){50}}

\put(202,52){$\bullet$}
\put(15,10){$L_+$}
\put(110,10){$L_{-}$}
\put(200,10){$L_{\times}$}

\end{picture}
\caption{The singular links $L_{+}$, $L_{-}$ and $L_{\times}$.}\label{skeinsfig}
\label{fig4}
\end{center}
\end{figure}

Let ${\mathcal L}_{\mathcal S}$ denote the set of oriented singular links.  Oriented singular links are represented by singular braids, which form the singular braid monoids  ${\mathcal S}B_n$ \cite{bae, bi, sm}.  ${\mathcal S}B_n$ is generated by the classical braiding generators $\sigma_i$ with their inverses and by the elementary singular braids $\tau_i$ which are not invertible.  In \cite{jula3} the following monoid homomorphism was constructed:
\begin{equation}\label{eta}
\begin{array}{cccl}
\eta : & {\mathcal S} B_n & \longrightarrow & {\rm Y}_{d,n}(u) \\
 & \sigma_i  & \mapsto & g_i \\
 & \tau_i    & \mapsto & p_i = e_i(1 + g_i)
\end{array}
\end{equation}
In view of the elements $p_i$ the quadratic relations (\ref{quadr}) may be rewritten as: $g_i^2 = 1 + (u-1)p_i$. Using now  the singular braid equivalence \cite{ge}, the map $\eta$ and the specialized Juyumaya trace ${\rm tr}_D$ we obtain isotopy invariants for oriented singular links \cite{jula3}:

\begin{thm}[\cite{jula3}]\label{sing} 
 Given a solution $X_D$ of the E--system, for any oriented singular link $\widehat{\alpha} \in {\mathcal L}_{\mathcal S}$, where $\alpha \in {\mathcal S}B_n$, we define:
$$
H_D (\widehat{\alpha}) = P_D^{n-1} (\sqrt{\lambda_D})^{\epsilon(\alpha)} \left({\rm tr}_D \circ \eta \right)(\alpha)
$$
where $P_D, \lambda_D$ are as defined in \S\ref{framedknots},  $\eta$ as defined in (\ref{eta}) and where $\epsilon(\alpha)$ is defined as follows: Let $\alpha = \mu_1^{\epsilon_1}\mu_2^{\epsilon_2}\ldots \mu_m^{\epsilon_m}$ where $\mu_j\in\{\sigma_i, \tau_i\,;\,1\leq i\leq n-1\}$. Then $\epsilon(\alpha) := \epsilon_1 +\ldots +\epsilon_m$, where $\epsilon_j = +1$ or $-1$ if  $\mu_j = \sigma_j$ and $\epsilon_j = +1$ if $\mu_j = \tau_j$.  Then the map  $H_D (u,\lambda_D)$ is a 2--variable isotopy invariant of oriented singular links.
\end{thm}
Moreover, in the image $\eta ({\mathcal S} B_n)$ the following relations hold:
\begin{equation}\label{gipi}
g_i^{-1} - g_i = (u^{-1}-1)p_i
\end{equation}
which give rise to the following skein relation:
\begin{equation}\label{skeins}
\sqrt{\lambda_D}\, H_D (L_{-})  -  \frac{1}{\sqrt{\lambda_D}}\, H_D (L_{+})
= \frac{u^ {-1}-1}{\sqrt{\lambda_D}}\, H_D (L_{\times})
\end{equation}
where  $L_{+}$, $L_{-}$ and $L_{\times}$ are diagrams of three oriented singular links, which are identical except near one crossing, where they are as depicted in Figure~\ref{skeinsfig}.

\smallbreak

For further details the reader is referred to \cite{jula3}. We note that there are some differences from \cite{jula3} in the signs in (\ref{eta}), (\ref{gipi}) and  (\ref{skeins}); this is due to the change of sign in the quadratic relation (\ref{quadr}).

\subsection{}\label{trans}

Another category of knots which seems to be naturally related to the Yokonuma--Hecke algebras is the category of {\it transverse knots}, for the following reasons: transverse knots are naturally  framed  and oriented, and their equivalence is more refined than framed knot isotopy.   In 1983 D. Bennequin \cite{Be} noted that the closed braid presentation of knots is convenient for describing transverse knots. Further, S.~Orevkov and V.~Shevchishin \cite{OrSh}  and independently N.~Wrinkle \cite{Wr} gave  a transverse analogue of the Markov theorem,
 comprising conjugation in the framed braid groups and {\it positive only}  stabilizations and destabilizations: $ \alpha \sim \alpha \,{\sigma}_n \in {\mathcal F}_{n+1}$, where $\alpha \in {\mathcal F}_n$.

\smallbreak
 Now, rule~(2) of the trace ${\rm tr}$ (Theorem~\ref{trace}) tells us that it respects positive stabilizations but property~(\ref{nofactor}) tells us that ${\rm tr}$ does not behave well under negative stabilizations. 
Using the above, an invariant $M (u,\lambda, x_1, \ldots , x_{d-1})$ for transverse  knots has been constructed in \cite{chjajukala}, which  coincides with the  invariant  $\Gamma_D (u,\lambda_D)$ of framed  knots when $(x_1, \ldots , x_{d-1})$ specializes to a solution  of the E--system. However, as it turns out, the  invariants $M (u,\lambda, x_1, \ldots , x_{d-1})$ are only topological invariants of  framed knots. For details see  \cite{chjajukala}.

\subsection{}\label{comp}

Due to the quadratic relation (\ref{quadr}) it is very difficult to do computations in the Yokonuma--Hecke algebras and the trace ${\rm tr}$. In http://www.math.ntua.gr/\~{}sofia/yokonuma/index.html the reader can find a program by Sergei Chmutov and Konstantinos Karvounis, together with instructions for using it.

\section{Framization of the Temperley--Lieb algebra}\label{tl}

In this section we define possible framizations for the Temperley--Lieb algebra and we indicate one of them as our favourite. The results of this section are contained  in \cite{gojukola} and \cite{gojukolaf}.

\subsection{}
The Temperley--Lieb  algebra can be defined  in several ways. Here, it is  convenient to define  the {\it Temperley--Lieb algebra} ${\rm TL}_n(u)$ as the quotient of the Hecke algebra ${\rm H}_n(u)$  by the two--sided ideal generated by the {\it Steinberg elements} \cite{jo}:
$$
h_{i,i+1} :=  h_ih_{i+1}h_i + h_{i+1}h_i + h_i h_{i+1} + h_i  + h_{i+1} + 1,
$$
that is:
\begin{equation}\label{tln}
{\rm TL}_n(u) = \frac{{\rm H}_n(u)}{\langle h_{i,i+1}  \,  ; \ \text{ all } i \rangle}. 
\end{equation}

V.F.R. Jones recovered his Markov trace on the Temperley--Lieb algebra from the Ocneanu trace $\tau$  \cite{jo} as follows: He required first that $\tau$ is zero on the Steinberg elements. This yielded two values for the parameter $\zeta$, namely:
\begin{equation}\label{jonval}
\zeta = - \frac{1}{u+1}, \ \text{ where } u\neq -1 \quad \mbox{and} \quad \zeta=-1.
\end{equation}
He then showed that for these values $\tau$ is zero on all elements of the defining ideal of ${\rm TL}_n(u)$.

\smallbreak
Finally, the Jones polynomial  $V(u)$ is obtained from the Homflypt polynomial (recall \S~\ref{homflypt}) by taking  
$$
\zeta = -\frac{1}{u+1}.
$$
Hence $\lambda = u$. That is, 
$$
V(u) = P(u,u).
$$   

In order to define a framization of the Temperley--Lieb algebra, the most natural approach is to take an appropriate quotient of the Yokonuma--Hecke algebra, in analogy with the classical case. So, we try to define in ${\rm Y}_{d,n}(u)$ elements analogous to the $h_{i,i+1}$'s as generators of a two--sided ideal. Now, having in mind that the multiplication rules of the Hecke algebra are governed by the symmetric group and those of ${\rm Y}_{d,n}(u)$  are governed by the group $C_{d, n}$ (recall Eq.~\ref{cdn}), we have two obvious possibilities for generating elements, which give rise to  two natural  candidates for the framization of the Temperley--Lieb algebra.

\subsection{}

For the first possibility we define for all $i$ the elements:
\begin{equation}\label{gij}
g_{i,i+1} : = g_ig_{i+1}g_i + g_{i+1}g_i + g_i g_{i+1} + g_i  + g_{i+1} + 1.
\end{equation}
Then, we have: 
\begin{defn}[\cite{gojukola}]\label{ytldef} \rm For $n\geq 3$, the  algebra  {\it Yokonuma--Temperley--Lieb}, denoted ${\rm YTL}_{d,n}(u)$, is defined as the quotient:
$$
{\rm YTL}_{d,n}(u) := \frac{{\rm Y}_{d,n}(u)}{\langle g_{i,i+1} \,  ; \ \text{ all }  i \rangle}. 
$$
\end{defn}

 In   \cite{chpou} the representations of ${\rm YTL}_{d,n}(u)$ are determined, as well as the dimension and an explicit basis. In \cite{gojukola} ${\rm YTL}_{d,n}(u)$ is studied as a knot algebra. A presentation with non--invertible generators is given. Also, the necessary and sufficient conditions are established for the Juyumaya trace ${\rm tr}$ on ${\rm Y}_{d,n}(u)$ to pass to the quotient algebra ${\rm YTL}_{d,n}(u)$. Indeed we have: 

\begin{thm}[Theorem~5 \cite{gojukola}]\label{ytlthmgen}
The trace ${\rm tr}$ passes to the quotient ${\rm YTL}_{d,n}(u)$ if and only if the $x_i$'s  are solutions of the E--system and one of the two cases holds: 
\begin{enumerate}
\item [(i)] 
For some  $0 \leq m_1 \leq d-1$ the $x_\ell$'s are expressed as:
\[
x_\ell=\exp(\ell m _1) \quad (0 \leq \ell \leq d-1).
\]
In this case  the $x_\ell$'s are $d^{th}$ roots of unity and \ 
$z=-\frac{1}{u+1}$ \ or \ $z=-1$.
\item [(ii)] 
For some $m_1, m_2$ such that $0 \leq m_1 \neq m_2 \leq d-1$ the $x_\ell$'s are expressed as: 
$$
x_\ell=\frac{1}{2}\left(\exp(\ell m _1) +\exp(\ell m _2) \right) \quad (0 \leq \ell \leq d-1).
$$
In this case 
 we have \ $z=-\frac{1}{2}$.
\end{enumerate}
\end{thm}

The cases where  $z=-1$ or $z=-1/2$ are of no topological interest.  The only interesting case is case $(i)$ for $z=-\frac{1}{u+1}$. In this case the $x_\ell$'s are $d^{th}$ roots of unity,  which is equivalent to $E=1$ and $\vert D \vert = 1$  \cite[Appendix]{jula5}.  By the discussion at the end of \S~\ref{classical} the invariants $\Delta_D(u,\lambda_D)$ coincide with the Homflypt polynomial.  Further, $z=-\frac{1}{u+1}$ implies $\lambda_D = u$.  Hence, the invariants  derived for classical, framed and singular  knots are $\Gamma_D(u,u)$, $\Delta_D(u,u)$ and  $H_D(u,u)$  respectively, recall \S \ref{invariants}. 
\begin{rem} \rm 
In particular, the invariants: 
$$
V_D(u) := \Delta_D(u,u)
$$ 
for classical knots all {\it coincide with the Jones polynomial}. From our point of view, this is the characteristic property of the algebra  ${\rm YTL}_{d,n}(u)$.
\end{rem} 

To recapitulate, the conditions on the trace ${\rm tr}$ are too strong, resulting in a trivial framization of the Jones polynomial. For this reason the algebra ${\rm YTL}_{d,n}(u)$ is discarded as a framization of ${\rm TL}_n(u)$.

\subsection{}

For the second possibility we define for all $i$ the following elements, where the framings are now also involved.
$$
c_{i,i+1} : = \sum_{\alpha,\beta,\gamma \in {\mathbb Z}/d{\mathbb Z}} t_i^{\alpha} t_{i+1}^{\beta} t_{i+2}^{\gamma} \, g_{i,i+1},
$$
where $g_{i,i+1}$ as defined in (\ref{gij}). Then, we obtain: 
\begin{defn}[\cite{gojukolaf}] \rm
The {\it Complex Reflection Temperley--Lieb algebra}, denoted ${\rm CTL}_{d,n}(u)$, is defined as the quotient: 
$$
{\rm CTL}_{d,n}(u) := \frac{{\rm Y}_{d,n}(u)}{\langle  c_{i,i+1} \,  ; \ \text{ all }  i  \rangle}. 
$$
\end{defn}
In \cite{gojukolaf}  the necessary and sufficient conditions are determined for the  trace ${\rm tr}$  to pass to ${\rm CTL}_{d,n}(u)$. More precisely, we have the following:
\begin{thm}[\cite{gojukolaf}]\label{ctlthm}
The trace ${\rm tr}$ passes to the quotient ${\rm CTL}_{d,n}(u)$ if and only if the parameter $z$ and the $x_i$'s are related through the equation:
 $$
 (u+1)z^2\sum_{k\in \mathbb{Z}/d\mathbb{Z}} x_k  +  (u+2)z\sum_{k\in \mathbb{Z}/d\mathbb{Z}} E^{(k)}+ \sum_{k \in \mathbb{Z}}{\rm tr}(e_1^{(k)} e_2)=0
 $$ 
where
$$
e_1^{(k)} := \frac{1}{d} \sum_{s=0}^{d-1} t_1^{k+s} t_{2}^{d-s} \quad \text{ and } \quad  E^{(k)} :={\rm tr}(e_1^{(k)}) \qquad (0\leq k \leq d-1)
$$
and where $e_1^{(0)} = e_1$ and $ E^{(0)} = E$.
\end{thm}

 \begin{rem} \rm Contrary to the case of ${\rm YTL}_{d,n}(u)$, the conditions of Theorem~\ref{ctlthm} are too relaxed  on the trace parameters $x_i$. 
However, in order to obtain framed, classical or singular link invariants from the algebras ${\rm CTL}_{d,n}(u)$ one has to impose the E--condition on the $x_i$'s. This in turn leads to the following values for $z$ (see \cite{gojukolaf}): 
$$
z= -\frac{1}{(u+1)|D|} \ \text{ or } \ z=-\frac{1}{|D|}.
$$
 Again, the interesting value is $z= - 1 / (u+1)|D|$ and, as it turns out, the link invariants one obtains  coincide with those from the algebras ${\rm FTL}_{d,n}(u)$ that we define next.
\end{rem}

\subsection{}

From the above constructions and from the knot algebras point of view it follows that it would be more sensible to consider an intermediate algebra between ${\rm CTL}_{d,n}(u)$ and ${\rm YTL}_{d,n}(u)$, for which the conditions for the trace ${\rm tr}$ to pass through include explicitely all solutions of the E--system.
 Indeed, for all $i$ we define the elements:
\begin{equation}\label{ftlideal}
r_{i,i+1} := \sum_{\substack{{\alpha,\beta,\gamma \in {\mathbb Z}/d{\mathbb Z}}\\ \alpha+\beta+\gamma=0 }} t_i^\alpha t_{i+1}^{\beta}  t_{i+2}^\gamma \, g_{i,i+1}  = e_i e_{i+1}\, g_{i,i+1}
\end{equation}
We now define:
\begin{defn}[\cite{gojukolaf}]  \rm
The {\it framization of the Temperley--Lieb algebra}, denoted ${\rm FTL}_{d,n}(u)$,  is defined as the quotient:
$$
{\rm FTL}_{d,n}(u) := \frac{{\rm Y}_{d,n}(u)}{\langle  r_{i, i+1}  \,  ; \ \text{ all }  i \rangle}. 
$$
\end{defn}
For the algebra ${\rm FTL}_{d,n}(u)$ we have determined in \cite{gojukolaf} the  necessary and sufficient conditions on the trace parameters  $z$ and the $x_i$'s for the trace ${\rm tr}$ to pass to ${\rm FTL}_n(u)$.  In order to  state our result we need to introduce the following notation: given a sequence $(x_1, \ldots, x_{d-1})$ of $d-1$ complex numbers, we denote by $x$ the function on ${\mathbb Z}/d{\mathbb Z}$ with values in ${\mathbb C}$, such that $0$ is mapped to $1$ and $i$ in mapped to $x_i$, for $1\leq i\leq d-1$. Further, we denote $\widehat{x}$ the Fourier transform of $x$ (cf. \cite[\S 3]{gojukolaf}). 

We then have the following:

\begin{thm}[\cite{gojukolaf}]\label{trpassesFTL}
The trace {\rm tr} passes to ${\rm FTL}_{d,n}(u)$ if and only if the trace parameters satisfy 
$$
x_m= -z\left(\sum_{k\in D_1}  {\rm exp}(km)
+
\sum_{k\in D_2}   {\rm exp}(km)\right) \qquad {\text and} \qquad z=  - \frac{1}{\vert D_1 \vert + (u+1)\vert D_2 \vert}
$$
where the  disjoint union $D_1\cup D_2$ is  the support of the Fourier transform $\hat{x}$ of $x$, 
$$
D_1: = \{ k\in {\mathbb Z}/d{\mathbb Z}\,;\, y_k = - dz\} \ ,
\quad
D_2:= \{ k\in {\mathbb Z}/d{\mathbb Z}\,;\, y_k = - dz(u+1)\} \ ,
$$
and  the $y_k$'s denote the values of  $\hat{x}$.
\end{thm}
In particular, the above theorem implies that whenever the trace parameters  are solutions of the E--system, then the trace ${\rm tr}$ passes to ${\rm FTL}_{d,n}(u)$.
Indeed, we have the following corollary.

\begin{cor}[\cite{gojukolaf}]\label{ftlpass}
In the case where  one of the sets $D_1$ or $D_2$ is  the empty set we obtain that the values in the above theorem  become  solutions  of the E--system.
More precisely,  in the case  where $D_1$ is  the empty set, the 
 $x_m$'s become  the solutions of the E--system parametrized by $D_2$ and $z$ takes the value  $z = -1/(u+1)\vert D_2\vert$. In the case where $D_2$ is  the empty set we have that the $x_m$'s become  the solutions of the E--system parametrized by $D_1$ and $z$ takes the value  $z= -1/ \vert D_1\vert$.
\end{cor}

The above Corollary allows one to define a non--trivial version of a framed Jones polynomial in the same manner as the Jones polynomial coincides with $P(u,u)$, the specialization   of the Homflypt polynomial  $P(u,\lambda)$ for $z=-1/(u+1)$. More precisely, we give the following definition, which is possible by Corollary~\ref{ftlpass} and Eq.~\ref{Econdition}.

\begin{defn}[\cite{gojukolaf}]\rm 
The {\it framed Jones polynomial} ${\mathcal V}_D(u)$ is defined as the specialization of  the polynomial $\Gamma_D(u, \lambda_D)$ at the value: 
$$
z = -\frac{1}{(u+1)\vert D\vert}. 
$$
Hence:
$$
{\mathcal V}_D(u):= \Gamma_D(u, u).
$$
\end{defn}

Similarly, one derives the invariants $\Delta_D(u,u)$ and  $H_D(u, u)$  for  classical and  singular links respectively from the invariants $\Delta_D(u,\lambda_D)$ and  $H_D(u,\lambda_D)$  (recall \S \ref{invariants}) by specializing at $z = -1/(u+1)\vert D\vert$. 

\begin{rem} \rm 
The invariants for classical links $\Delta_D(u,u)$  may prove to be topologically equivalent with the Jones polynomial, in analogy with $\Delta_D(u,\lambda_D)$ and the  Homflypt polynomial (recall discussion at the end of \S~\ref{classical}). 
In case they are not, one could compare corresponding 3-manifold invariants with the Witten invariants from the Jones polynomial, see \cite{we}.
\end{rem}

\subsection{}
All three quotient algebras ${\rm YTL}_{d,n}(u)$, ${\rm FTL}_{d,n}(u)$ and ${\rm CTL}_{d,n}(u)$ equipped with the Markov traces are interesting on their own right. They are related via the  following algebra epimorphisms:
$$
{\rm Y}_{d,n}(u) \twoheadrightarrow {\rm CTL}_{d,n}(u) \twoheadrightarrow {\rm FTL}_{d,n}(u) \twoheadrightarrow {\rm YTL}_{d,n}(u),
$$
which follow from inclusions of their defining ideals \cite{gojukolaf}.

\section{Framization of the Hecke--related algebras of type $B$}\label{b}

\subsection{} 

The {\it Artin group of type $B$}, denoted $B_{1,n}$, is related to the following Dynkin diagram:

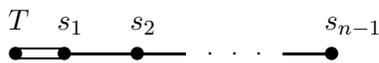
\begin{figure}[ht!]
\begin{center}
{\tt    \setlength{\unitlength}{0.92pt}
\begin{picture}(300,50)
\put( 80, 20){\circle*{5}}
\put( 80, 18){\line(1,0){20}}
\put( 80, 22){\line(1,0){20}}
\put( 100, 20){\circle*{5}}
\put( 100, 20){\line(1,0){30}}
\put( 130, 20){\circle*{5}}
\put( 130, 20){\line(1,0){20}}
\put(160, 20){\circle*{1}}
\put(170, 20){\circle*{1}}
\put(180, 20){\circle*{1}}
\put(190, 20){\line(1,0){20}}
\put(210, 20){\circle*{5}}
\put( 77, 30){$T$}
\put( 97, 30){$s_1$}
\put( 127, 30){$s_2$}
\put(207, 30){$s_{n-1}$}
\end{picture}}
\caption{The Dynkin diagram of type $B$.}
\label{}
\end{center}
\end{figure}

\begin{figure}[ht!]
\begin{center}
\includegraphics[width=7cm]{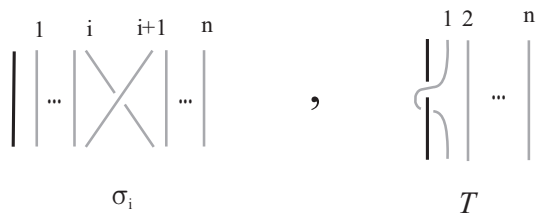}
\caption{The braiding generators and the loop generator of $B_{1,n}$.}
\label{braidgenrs}
\end{center}
\end{figure}

That is, $B_{1,n}$ is presented by the braiding generators $\sigma_1 ,\ldots, \sigma_{n-1}$ and the loop generator $T$ (see Figure~\ref{braidgenrs}), satisfying the braid relations (\ref{braidrels}) and the relations:
 $$
 \begin{array}{rcll}
 T \sigma_1 T \sigma_1 & = & \sigma_1 T \sigma_1 T &   \\
T \sigma_i & = & \sigma_i T & \text{ for} \quad 2\leq i \leq n-1.
\end{array}
 $$ 

\smallbreak

\begin{figure}
\begin{center}
\includegraphics[width=8cm]{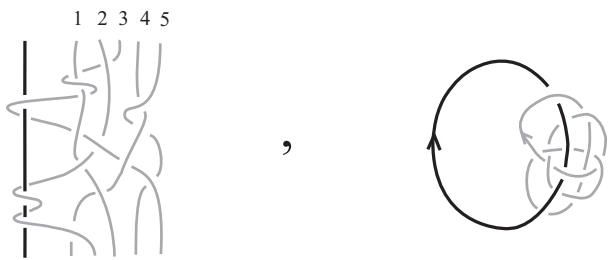}
\caption{A braid in $B_{1,n}$ and a link in the solid torus.}
\label{braidknot}
\end{center}
\end{figure}

Geometrically, a braid in $B_{1,n}$ has $n+1$ strands, with the first strand identically fixed  and the other $n$ strands numbered from 1 to $n$. See Figure~\ref{braidknot} for an example. Its closure is an oriented link in the solid torus, where the complement solid torus is represented by the closure of the fixed strand \cite{la1, la2}.
  
\begin{rem}\label{cylinderbraid} \rm
The Artin braid group of type $B$, $B_{1,n}$, is isomorphic to the affine Artin braid group of type $A$, denoted $\widetilde{B_n}$. Thus, another geometric interpretation for elements in  $B_{1,n}$ is as cylinder braids, that is, as braids in a thickened cylinder. Then, depending on how closure is defined (by simple closed arcs in the thickened cylinder or by endpoints identifications), they could give rise to oriented links in the solid torus or in the thickened torus.
\end{rem}
	
\smallbreak

\subsection{}\label{btypealgebras}

For $u,U \in {\mathbb C}\backslash \{0\}$, the classical {\it Iwahori--Hecke algebra of type $B$} \cite{DJ}, denoted here ${\rm H}_{1,n}(u,U)$, can be viewed as the  quotient of the  group algebra $\mathbb C B_{1,n}$  over the expressions: 
$$
\sigma^2_i  - (u-1)\sigma_i - u \qquad \text{ and } \qquad  T^2  - (U-1)T - U.
$$  
  
Further, for $u,u_1, \ldots, u_r \in {\mathbb C}\backslash \{0\}$, the {\it cyclotomic Hecke algebra of type $B$ and  degree $r$} \cite{arko,brma}, denoted here ${\rm H}_{1,n}(u,u_1, \ldots, u_r)$, can be defined  as the quotient of the group algebra  $\mathbb C B_{1,n}$  over the expressions: 
$$
\sigma^2_i  - (u-1)\sigma_i - u \qquad \text{ and } \qquad  (T -u_1)(T -u_2)\ldots(T -u_r).
$$  
  For $r=2$ the algebra ${\rm H}_{1,n}(u,u_1,u_2)$ can be proved to be isomorphic to ${\rm H}_{1,n}(u,U)$.
	
	\smallbreak
Finally, the {\it generalized Hecke algebra of type $B$} \cite{la2}, denoted ${\rm H}_{1,n}(u)$, is defined as the quotient of the group algebra  $\mathbb C B_{1,n}$  over the  expressions: 
$$
\sigma^2_i  - (u-1)\sigma_i - u.
$$    
The algebra ${\rm H}_{1,n}(u)$ was observed by T. tom Dieck \cite[Remark~1]{la2} to be isomorphic to the {\it affine Hecke algebra of type $A$}.

\smallbreak
In \cite{la1,gela} and \cite{la2} Markov traces have been constructed on all these algebras, giving rise to all possible analogues of the Homflypt polynomial for oriented links in the solid torus. The two rules of these traces, namely conjugation and Markov property, are analogous to the two rules of the Ocneanu trace (recall Theorem~\ref{ocneanu}). Then, there is also a third inductive rule that takes care of the loopings, namely:
\begin{equation} \label{btracerule}
\tau\big(a \, (h_n \ldots h_1 T^k h_1^{-1} \ldots h_n^{-1}) \big) = s_k \, \tau(a),
\end{equation}
where $\sigma_i$ corresponds to the braiding generator $h_i$ and $T$ corresponds to $T$, and where $a \in {\rm H}_{1,n}(u,U)$ or ${\rm H}_{1,n}(u,u_1, \ldots, u_r)$ or ${\rm H}_{1,n}(u)$, depending on the algebra we are in each time. For details see \cite{la1,gela} and \cite{la2}.

\subsection{}

From the above it is natural to try to define framizations of the Hecke algebra of type $B$, of the cyclotomic Hecke algebra of type $B$ and of the generalized Hecke algebra of type $B$, all $B$--type analogues of the Yokonuma--Hecke algebra, with the prospective  to obtain invariants of  links in the solid torus, framed, classical or singular,  analogous to those defined  by the Yokonuma--Hecke algebra. For this we define first:

\begin{defn}[Definition 4 \cite{jula6}] \rm
The {\it framed braid group of type} $B$, ${\mathcal F}_{1,n}$, is the group  presented by the generators $T, \sigma_1, \ldots, \sigma_{n-1}$ of $B_{1,n}$ together with the framing generators $t_1, \ldots , t_n$, subject to the relations  of $B_{1,n}$ together with the framing relations (\ref{framerels}) and the  relations:
$$
\begin{array}{rcl}
t_i T & = & T t_i \quad \text{for all}\quad  1\leq i\leq n.
\end{array}
$$
Geometrically, elements in ${\mathcal F}_{1,n}$ have framings on the $n$ numbered strands and upon closure they represent framed links in the solid torus.  The {\it $d$--modular framed braid group of type} $B$, denoted ${\mathcal F}_{d,1,n}$, is defined by adding to the above presentation of ${\mathcal F}_{1,n}$ the relations:
$$
 \begin{array}{rcll}
t_i^d & = & 1 & \text{ for}  \quad  1\leq i \leq n.
\end{array}
 $$
\end{defn}

\begin{rem}\label{cylinderframedbraid} \rm
By Remark~\ref{cylinderbraid} ${\mathcal F}_{1,n}$ is isomorphic to the affine framed braid group of type $A$, denoted $\widetilde{{\mathcal F}_n}$. So, elements in  ${\mathcal F}_{1,n}$ can be also interpreted as framed braids in a thickened cylinder. Then, upon different types of  closures, they  give rise to oriented framed links in the solid torus or in the thickened torus.
\end{rem}

\subsection{}

Define now the following elements in ${\mathbb C} {\mathcal F}_{d,1,n}$: 
$$
e_i:= \frac{1}{d}\sum_{m=0}^{d-1} t_i^{m}t_{i+1}^{d-m} \quad 
\text{for all}\quad 1\leq i\leq n-1.
$$
It is clear that the elements $e_i$ are idempotents.  The  $e_i$'s are represented geometrically as in Figure~\ref{ed1}, but with the addition in all terms of a first fixed strand with no framing.

\smallbreak

We shall now define framizations of the $B$--type algebras above. These definitions  have  all been given in  \cite[Definition~5]{jula6}.

\begin{defn}[\cite{jula6}] \rm
For $u, U \in {\mathbb C}\backslash \{0\}$  we define the {\it Yokonuma--Hecke algebra of  type $B$}, ${\rm Y}_{d,1,n}(u,U)$, as the algebra presented by the braiding generators $g_1, \ldots, g_{n-1}$, the loop generator $T$ and the  framing generators $t_1, \ldots, t_n$, subject to the relations (\ref{yh}) for the braiding and framing generators, together with the following extra relations:
\begin{equation}\label{brels}
 \begin{array}{rcll}
 T g_1 T g_1 & = & g_1 T g_1 T &   \\
T g_i & = & g_i T  & \text{ for} \quad 2\leq i \leq n-1 \\ 
T t_j  & = & t_j T  & \text{ for} \quad 1\leq j \leq n 
\end{array}
\end{equation}
and the quadratic relations: 
\begin{equation}\label{bquadr1}
g_i^2 = 1 + (u-1)e_i + (u-1)e_ig_i   \quad (1\leq i \leq n-1) 
\end{equation}
and
\begin{equation}\label{bquadr2}
T^2 = (U-1)T +  U .
\end{equation}
\end{defn}

The generators $g_i$ and $T$ are are easily seen to be invertible:
\begin{equation}\label{binvrs}
g_i^{-1} = g_i + (u^{-1} - 1)\, e_i + (u^{-1} - 1)\, e_i \, g_i \qquad \text{and} \qquad T^{-1} = U^{-1} T + (U^{-1} - 1).
\end{equation}
From the above presentation it is clear that the algebra ${\rm Y}_{d,1,n}(u,U)$ is a quotient of the modular framed braid group algebra ${\mathbb C} {\mathcal F}_{d,1,n}$ under the quadratic relations (\ref{bquadr1}) and (\ref{bquadr2}). This observation leads to diagrammatic interpretations for the elements of ${\rm Y}_{d,1,n}(u,U)$. 

\begin{defn}[\cite{jula6}] \rm
For $u, u_1, \ldots, u_r \in {\mathbb C}\backslash \{0\}$ we define the {\it cyclotomic Yokonuma--Hecke algebra of  type $B$ and of degree $r$}, ${\rm Y}_{d,1,n}(u, u_1, \ldots, u_r)$, as the algebra presented by the braiding generators $g_1, \ldots, g_{n-1}$, the loop generator $T$ and the  framing generators $t_1, \ldots, t_n$ subject to the relations (\ref{yh}) for the braiding and framing generators,  relations (\ref{brels}) for the loop generator, the quadratic relations (\ref{bquadr1}) for the braiding generators and the following polynomial relation for the loop generator $T$ in place of (\ref{bquadr2}): 
\begin{equation} \label{cyclotomic}
(T -u_1)(T -u_2)\ldots(T -u_r) = 0.
\end{equation}
\end{defn}

 Note  that the generators $g_i$ and $T$ are invertible. Also, that the algebra ${\rm Y}_{d,1,n}(u, u_1, \ldots, u_r)$ is clearly a quotient of the modular framed braid group algebra ${\mathbb C} {\mathcal F}_{d,1,n}$ under the quadratic relations (\ref{bquadr1}) and relation (\ref{cyclotomic}). 

\begin{defn}[\cite{jula6}] \rm
For $u \in  {\mathbb C}\backslash \{0\}$ we define the {\it generalized Yokonuma--Hecke algebra of  type $B$}, ${\rm Y}_{d,1,n}(u)$, as the algebra presented by the braiding generators $g_1, \ldots, g_{n-1}$, the loop generator $T$ and the  framing generators $t_1, \ldots, t_n$ subject to the relations (\ref{yh}) for the braiding and framing generators,   relations (\ref{brels}) for the loop generator and the quadratic relations (\ref{bquadr1}) for the braiding generators.
\end{defn}

\begin{rem} \rm
By \cite[Remark~1]{la2} and Remark~\ref{cylinderframedbraid}, the algebra ${\rm Y}_{d,1,n}(u)$ can be considered isomorphically as the {\it affine  Yokonuma--Hecke algebra of  type $A$}. See also \cite{chpa}.
\end{rem}

Note that for $d=1$ the  algebras ${\rm Y}_{1,1,n}(u,U)$, ${\rm Y}_{1,1,n}(u, u_1, \ldots, u_r)$  and ${\rm Y}_{1,1,n}(u)$ coincide with the corresponding algebras of $B$--type defined above. 
Note also that the algebras ${\rm Y}_{d,1,n}(u,U)$  and ${\rm Y}_{d,1,n}(u, u_1, \ldots, u_r)$ can be clearly viewed as quotients of ${\rm Y}_{d,1,n}(u)$ by the relations (\ref{bquadr2}) and (\ref{cyclotomic}) respectively. Namely, we have algebra epimorphisms:
$$
{\mathbb C} {\mathcal F}_{d,1,n} \twoheadrightarrow {\rm Y}_{d,1,n}(u) \twoheadrightarrow {\rm Y}_{d,1,n}(u, u_1, \ldots, u_r).
$$

On all these algebras unique Markov traces can be constructed with the three rules analogous to the rules of the Markov traces discussed in \S~\ref{btypealgebras}, together with a fourth inductive rule that takes care of the framings, analogous to rule~(3) of the Juyumaya trace ${\rm tr}$ (recall Theorem~\ref{trace}). Such traces are constructed in \cite{chpa2}. Further, invariants of framed links, classical links and singular links in the solid torus can be constructed, up to necessary conditions analogous to the E--condition. These conditions are given  in~\cite{chpa2}, where the representation theory of all these algebras is also studied.

\section{Framization of the BMW algebra}\label{bmw}

\subsection{}

Let $l,m \in  {\mathbb C}\backslash \{0\}$.  Birman  and Wenzl \cite{biwe} and  simultaneously but independently Murakami \cite{mu}  defined  a unital associative algebra $C_n=C_n(l,m)$ which is known as the {\it Birman--Wenzl--Murakami algebra} or simply the {\it BMW algebra}. The algebra $ C_n$  is defined by  two sets of generators: the `braiding' generators $g_1, \ldots , g_{n-1}$ and the `tangle' generators $h_1, \ldots , h_{n-1}$, satisfying:  
 the  braid relations (\ref{braidrels}) among the $g_i$'s, together with the relations:
\begin{equation}\label{bmw1}
 \begin{array}{rcll}
 g_ih_i & = & l^{-1}h_i &  \text{ for} \quad 1 \leq i \leq n-1  \\
h_ig_{i\pm 1} h_i  & = & lh_i  & \text{ for} \quad \text{all } i  \\
g_j h_i & = & h_i g_j &  \text{ for} \quad \vert i-j \vert \geq 2 
\end{array}
\end{equation}
and the quadratic relations:
\begin{equation}\label{bmwquadr}
 g_i^2   =  1 -mg_i + ml^{-1}h_i   \qquad (1\leq i \leq n-1).
\end{equation}
For diagrammatic interpretations of the `tangle' generators $h_i$ the reader may view Figure~\ref{fbmw1}, where the framings 0 and 1 should be ignored and where $h_i$ corresponds to the two horizontal arcs joining the $i$ and $i+1$ endpoints.

From the defining relations of $C_n$ we deduce that the $g_i$'s are invertible:
\begin{equation}\label{bmwinvrs}
 g_i^{-1}  = g_i- m h_i + m
\end{equation}
 and also the following important relations:
\begin{equation}\label{bmw2}
\begin{array}{rcll}
h_ig_i & = & l^{-1}h_i &  \text{ for} \quad 1 \leq i \leq n-1  \\
h_ih_j  & = & h_jh_i & \text{ for} \quad \vert i-j\vert \geq 2  \\ 
h_i^2  & = &y h_i  & \text{ for} \quad 1\leq i \leq n-1
\end{array}
\end{equation}
where
$$
y:= 1 +\frac{l^{-1}-l}{m}.
$$

The algebra $C_n$  is a quotient of the classical braid group algebra ${\mathbb C} B_n$. To see this consider  $C_n$  generated by the $g_i$'s only, and view Eq.~\ref{bmwinvrs} as the defining relations for the $h_i$'s.  Further, the element $h_i$ can be seen represented in the category of $(n,n)$--tangles  as the elementary tangle consisting in two curved parallel horizontal arcs joining the endpoints $i$ and $i+1$ at  the top and at the bottom of the otherwise identity tangle (see Figure~\ref{fbmw1}). 
The algebra  $C_n$  is related to the {\it Kauffman polynomial} invariant for  classical knots \cite{ka}.

\subsection{}
In \cite{jula6} the framization of the BMW algebra has been introduced. We shall decribe it here briefly.

\begin{defn}\rm
Let $y_0:= y$ and $y_1, \ldots , y_{d-1}$ in ${\mathbb C}\backslash \{0,1\}$.  The {\it framization of the BMW algebra $C_n$}, denoted $F_{d,n}= F_{d,n}(l,m, y_0,\ldots ,y_{d-1})$,  is defined through three sets of generators: the two sets of generators of the algebra $ C_n$ given above, together with the framing generators   $t_1, \ldots ,t_n$, satisfying all  defining relations of $ C_n$ except for the quadratic relations (\ref{bmwquadr}), which are replaced by the following quadratic relations:
\begin{equation}\label{fbmwquadr}
g_i^2 = (1 -m) - m \, e_i(g_i-1)  + m \, l^{-1}h_i \qquad  (1 \leq i \leq n-1)
\end{equation}
where $e_i$ as defined in (\ref{ei}), 
together with the following relations for the framing generators:
\begin{equation}\label{fbmwrels}
\begin{array}{rcll}
t_i^d = 1 & \text{ and} &  t_i t_j= t_j t_i &  \text{ for} \quad 1 \leq i, j \leq n  \\
t_jg_i  & = &    g_it_{s_i(j)} &  \text{ for} \quad1\leq i \leq n-1 \ \ \& \ \ 1\leq j \leq n \\
t_i h_i  & = & t_{i+1} h_i & \text{ for} \quad 1 \leq i \leq n-1  \\ 
 h_i t_i  & = & h_i t_{i+1}  & \text{ for} \quad 1\leq i \leq n-1  \\
h_i t_i^k h_i  & = & y_k \, h_i   & \text{ for} \quad 1\leq i \leq n-1 \ \& \ 0\leq k\leq d-1  \\
h_i t_j  & = & t_j h_i  & \text{ for} \quad j\neq i, i+1
\end{array}
\end{equation}
 where $s_i(j)$ is the effect  of  the transposition $s_i=(i, i+1)$ on $j$.
\end{defn}
Note  that for $d=1$ we have $e_i=1$, hence $F_{1,n}$ coincides with $C_n$. Also, the elements $g_i$ are invertible \cite[Proposition~1]{jula6}:
\begin{equation}\label{fbmwinvrs}
g_i^{-1} = \frac{1}{1-m}g_i - \frac{m}{1-m}g_ie_i - mh_i + me_i.
\end{equation}

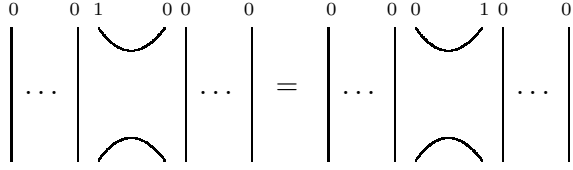
\begin{figure}
\begin{center}
\begin{picture}(220,60)
\put(-1,55){\tiny{$0$}}
\put(22,55){\tiny{$0$}}
\put(31,55){\tiny{$1$}}
\put(57,55){\tiny{$0$}}
\put(64,55){\tiny{$0$}}
\put(88,55){\tiny{$0$}}
\qbezier(0,0)(0,30)(0,50)
\put(5,25){$\ldots$}
\qbezier(25,0)(25,30)(25,50)
\qbezier(33 ,50)(45,33)(58,50)
\qbezier(33 ,0)(45,17)(58,0)
\qbezier(66,0)(66,30)(66,50)
\put(71,25){$\ldots$}
\qbezier(91,0)(91,30)(91,50)
\put(100, 25){$=$}
\put(119,55){\tiny{$0$}}
\put(142,55){\tiny{$0$}}
\put(151,55){\tiny{$0$}}
\put(177,55){\tiny{$1$}}
\put(184,55){\tiny{$0$}}
\put(208,55){\tiny{$0$}}
\qbezier(120,0)(120,30)(120,50)
\put(125,25){$\ldots$}
\qbezier(145,0)(145,30)(145,50)
\qbezier(153 ,50)(165,33)(178,50)
\qbezier(153 ,0)(165,17)(178,0)
\qbezier(186,0)(186,30)(186,50)
\put(191,25){$\ldots$}
\qbezier(211,0)(211,30)(211,50)
\end{picture}
\caption{The relation $t_ih_i= t_{i+1}h_i$.}
\label{fbmw1}
\end{center}
\end{figure}

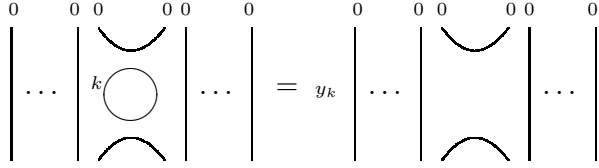
\begin{figure}
\begin{center}
\begin{picture}(220,60)
\put(-1,55){\tiny{$0$}}
\put(22,55){\tiny{$0$}}
\put(31,55){\tiny{$0$}}
\put(57,55){\tiny{$0$}}
\put(64,55){\tiny{$0$}}
\put(88,55){\tiny{$0$}}
\qbezier(0,0)(0,30)(0,50)
\put(5,25){$\ldots$}
\qbezier(25,0)(25,30)(25,50)
\put(45,25){\circle{20}}
\qbezier(33 ,50)(45,33)(58,50)
\qbezier(33 ,0)(45,17)(58,0)
\qbezier(66,0)(66,30)(66,50)
\put(71,25){$\ldots$}
\put(30,27){\tiny{$k$}}
\qbezier(91,0)(91,30)(91,50)
\put(100, 25){$=$}
\put(115,25){\tiny{$y_k$}}
\put(129,55){\tiny{$0$}}
\put(152,55){\tiny{$0$}}
\put(161,55){\tiny{$0$}}
\put(187,55){\tiny{$0$}}
\put(194,55){\tiny{$0$}}
\put(218,55){\tiny{$0$}}
\qbezier(130,0)(130,30)(130,50)
\put(135,25){$\ldots$}
\qbezier(155,0)(155,30)(155,50)
\qbezier(163 ,50)(175,33)(188,50)
\qbezier(163 ,0)(175,17)(188,0)
\qbezier(196,0)(196,30)(196,50)
\put(201,25){$\ldots$}
\qbezier(221,0)(221,30)(221,50)
\end{picture}
\caption{The relation $h_it_i^kh_i= y_kh_i$.}
\label{fbmw2}
\end{center}
\end{figure}

The algebra $F_{d,n}$ can be viewed as a quotient of the modular framed braid group algebra ${\mathbb C} {\mathcal F}_{d,n}$. To see this, in analogy to the classical case  we exempt the $h_i$'s from the set of generators for the algebra $F_{d,n}$ and  we consider Eq.~\ref{fbmwinvrs}  as the defining relations for the $h_i$'s.
So, elements in the algebra $F_{d,n}$ can be viewed as framed $(n,n)$--tangles, with framings modulo $d$. In this context, Figures~\ref{fbmw1} and~\ref{fbmw2} illustrate two of the relations of (\ref{fbmwrels}).

\smallbreak
In \cite{jula6} it is shown that $F_{d,n}$ is finite dimensional. We also have the following important result.

\begin{prop}[Proposition~2 \cite{jula6}]\label{quintic}
The elements $g_i$  satisfy  the quartic relation:
$$
 g_i^4 +mg_i^3 +(m-2)g_i^2+m (m-1)g_i - (m-1) = ml^{-1}\left(m+l^{-2}-1\right)h_i
$$
and this is of minimal degree not containing the framing generators $t_i$. Also, they satisfy the
 \lq closed\rq \ quintic relation:
$$
\left(x-l^{-1}\right)\big(x^4 +mx^3 + (m-2)x^2 +m (m-1)x -(m-1)\big)= 0
$$
and this is of minimal degree not containing the generators $t_i$ and $h_i$.
Moreover:
$$
x^4 +mx^3 + (m-2)x^2 +m (m-1)x -(m-1)= (x^2+mx-1)(x^2+ m-1).
$$
\end{prop}
Finally we have the following result:

\begin{prop}\cite[Proposition~3]{jula6}
Any element in $F_{d,n}$ can be written as a  ${\mathbb C}$--linear combination of monomials of the form $\alpha f\beta$, where $\alpha$ and $\beta$ are monomials in $1, g_1,\ldots, g_{n-2}$, $h_1, \ldots , h_{n-2}$, $t_1, \ldots ,t_{n-1}$ and 
$f\in X_n:=\{t_n^s, \ g_{n-1}, \ t_{n-1}^sh_{n-1}t_{n-1}^r \quad \vert \quad 0 \leq r,s \leq d-1 \}$.
\end{prop}
This result is in the right direction for defining Markov traces on $F_{d,n}$  via inductive rules.

\section {The framization of the singular Hecke algebra}\label{sin} 

A definition of the {\it singular Hecke algebra}, denoted $S{\rm H}_n(u)$, was proposed by Paris and Rabenda \cite{para}. This algebra is a  complex associative unital algebra defined  as the quotient of the  algebra ${\mathbb C}SB_n$ of the singular  braid monoid  $SB_n$ over the Hecke algebra quadratic relations (\ref{hecke}). 

Recall that  $SB_n$ was introduced in different contexts by Baez \cite{bae}, Birman \cite{bi} and Smolin \cite{sm} and it is defined by the classical elementary braids $\sigma_i$ with their inverses $\sigma_i^{-1}$, $i=1,\ldots, n-1$, which are subject to the braid relations (\ref{braidrels}), and by the elementary singular braids $\tau_i$, $i=1,\ldots, n-1$, together with the following relations:
\begin{equation}\label{hesin}
\begin{array}{rccl}
\lbrack\sigma_i, \tau_j\rbrack
=
 \lbrack\tau_i, \tau_j\rbrack & =  & 0 & \text{for}\quad \vert i-j\vert  >1 \\
   \lbrack\sigma_i, \tau_i \rbrack
& = &
0  & \text{for} \quad \text{all } i \\
 \sigma_i \sigma_j \tau_i
 &  = &
   \tau_j \sigma_i \sigma_j  & \text{for}\quad \vert i-j \vert = 1
\end{array}
\end{equation}

Paris and Rabenda \cite{para} constructed also a universal Markov trace on these algebras, which lead to an  invariant for singular knots, which is basically  equivalent to the singular link invariant of Kauffman and Vogel defined in \cite{kavo}.

\smallbreak
We can now proceed with the following definition.

\begin{defn}\rm
The {\it framization of the algebra $S{\rm H}_n(u)$}, denoted $FS_{d,n}(u)$, is defined as  the unital  associative algebra over $\mathbb C$, defined through three sets of generators: $g_1, \ldots, g_{n-1}$,  $\tau_1, \ldots, \tau_{n-1}$ corresponding to the two sets of generators of the algebra $ S{\rm H}_n(u)$ given above, together with the framing generators $t_1, \ldots ,t_n$, satisfying all the Yokonuma--Hecke algebra relations, (\ref{yh}) and (\ref{quadr}), together with the above relations (\ref{hesin}) of $ SB_n(u)$, whereby $\sigma_i$ corresponds to $g_i$ and $\tau_i$ to $\tau_i$.
\end{defn}

For the algebra $FS_{d,n}(u)$ one needs to find appropriate inductive basis and define on it a Markov trace analogue to the one by Paris and Rabenda.

\section{Concluding note}

We presented  framizations of several knot algebras, starting from the example of the classical Iwahori--Hecke algebra, and discussed questions that need to be further investigated.  There are many more other knot algebras. For example, other quotients of the classical braid group, quotients of the virtual braid group \cite{ka2,kala}, or the Rook algebra \cite{bry}, which is related to the Alexander polynomial.

\end{document}